%% file: ms.tex
\pdfoutput=1
\documentclass{article}
\usepackage[utf8]{inputenc}
\usepackage[T1]{fontenc}
\usepackage{graphicx}
\usepackage{grffile}
\usepackage{tikz-cd}
\usepackage{longtable}
\usepackage{wrapfig}
\usepackage{rotating}
\usepackage[normalem]{ulem}
\usepackage{amsmath}
\usepackage{textcomp}
\usepackage{bm}
\usepackage{siunitx}
\usepackage{amssymb}
\usepackage{amsmath}
\usepackage{capt-of}
\usepackage{hyperref}
\usepackage{subcaption}
\usepackage{cancel}
\usepackage{amsthm}
\usepackage{color}

\usepackage{indentfirst}

\author{Jack Coughlin\footnote{Department of Applied Mathematics, University of Washington, Seattle, WA 98195, USA (johnbc@uw.edu).} \ \ and \ Jingwei Hu\footnote{Department of Applied Mathematics, University of Washington, Seattle, WA 98195, USA (hujw@uw.edu).}}
\date{\today}
\title{Efficient dynamical low-rank approximation for the Vlasov-Amp\`{e}re-Fokker-Planck system\footnote{This work was partially supported by NSF CAREER grant DMS-2153208, NSF CDS\&E grant CBET-1854829, and AFOSR grant FA9550-21-1-0358.}}

\newcommand{\rd}{\,\mathrm{d}} 

\newcommand{\bO}{\mathcal{O}}
\newcommand{\bX}{\mathbf{X}} 
\newcommand{\bS}{\mathbf{S}} 
\newcommand{\bV}{\mathbf{V}} 
\newcommand{\bK}{\mathbf{K}} 
\newcommand{\bD}{\mathbf{D}}

\newcommand{\Mone}{\mathcal{M}_1}
\newcommand{\Mtwo}{\mathbf{M}_2}
\newcommand{\Mthree}{\overline{\overline{\mathsf{M}}}_3}

\numberwithin{equation}{section}

\usepackage[numbers]{natbib}

\begin{document}

\maketitle

\begin{abstract}
  Kinetic equations are difficult to solve numerically due to their high
  dimensionality. A promising approach for reducing computational cost is the
  dynamical low-rank algorithm, which decouples the dimensions of the phase space by
  proposing an ansatz as the sum of separable (rank-1) functions in position and velocity respectively. The fluid
  asymptotic limit of collisional kinetic equations, obtained in the
  small-Knudsen number limit, admits a low-rank representation when written as
  $f = Mg$, where $M$ is the local Maxwellian, and $g$ is low-rank. We apply
  this decomposition to the Vlasov-Amp\`{e}re-Fokker-Planck
  equation of plasma dynamics, considering the asymptotic limit of strong
  collisions and electric field. We implement our proposed algorithm and
  demonstrate the expected improvement in computation time by comparison to an
  implementation that evolves the full solution tensor $f$. We also demonstrate
  that our algorithm can capture dynamics in both the kinetic regime, and in the
  fluid regime with relatively lower computational effort, thereby efficiently
  capturing the asymptotic fluid limit.
\end{abstract}

{\small {\bf Key words.}  dynamical low-rank integrator, Vlasov-Amp\`{e}re-Fokker-Planck model, high-field limit,
convolution, implicit-explicit scheme}

\section{Introduction}
Magnetohydrodynamics and multi-fluid systems of equations provide reasonable
descriptions of plasma dynamics across a wide range of parameter regimes.
However, in situations where one or more of the particle species' phase space
distributions is far from a Maxwellian, fluid models can fail to capture
relevant physics. The Vlasov equation, when coupled with Maxwell's equations of
electrodynamics, provides a more complete description of plasma dynamics in
these regimes \cite{goldstonIntroductionPlasmaPhysics1995}. However, the
numerical solution of kinetic models is quite costly in 2 or 3 dimensions, since
they are posed over 4 or 6 phase space dimensions, respectively. This prompts
the search for computational algorithms which can accelerate the solution of
kinetic equations.

One promising approach for accelerated kinetic algorithms comes from the
recognition that there is low-rank structure in certain solutions
of kinetic equations. By low-rank structure, we mean that an approximation of the
following sort can be successful:
\begin{equation}
  \label{eqn:lowrank_ansatz}
  f(x, v, t) \approx \sum_{ij}^{r} X_{i}(x, t) S_{ij}(t) V_{j}(v, t).
\end{equation}
Such an approximation will be useful only if the numerical rank, $r$, is small
compared to the number of degrees of freedom $N_{x}$ and $N_{v}$. The
inspiration for this form of approximation, called a low-rank approximation,
comes from linear algebra and the need to deal effectively with extremely large
data matrices. The idea is to capture most of the action of a data matrix with a
low-rank approximation which can require vastly less storage.
Bounds on the quality of the approximation are available in various norms.
Perhaps the best-known approximation of this kind is the truncated Singular
Value Decomposition, which is known to provide the best rank $r$ approximation
to a given matrix in the spectral norm \cite{trefethenNumericalLinearAlgebra1997a}.

Our setting is time-dependent kinetic equations, so it is not enough to be able
to compress a given phase space distribution. One must also be able to evolve
the distribution $f$ in an approximate form. This is made possible by the theory
of \emph{dynamical} low-rank approximation, which has been studied in the matrix
and tensor contexts in \cite{kochDynamicalLowRank2007},
\cite{kochDynamicalTensorApproximation2010}.
This method advances the representation \ref{eqn:lowrank_ansatz} by updating the
the bases $X_{i}, V_{j},$ and the matrix of singular values $S_{ij}$ at each
timestep. The point is to never form the full product of all three factors but to
evolve the factorized form directly.
A crucial innovation in this field that we make use of
is the projector-splitting integrator of
\cite{lubichProjectorsplittingIntegratorDynamical2014}. This integrator enables
a robust dynamical low-rank method which is insensitive to
``overapproximation'', or vanishingly small singular values in the
approximation. Another such integrator with comparable robustness properties is 
the recently proposed ``unconventional'' dynamical low-rank integrator of
\cite{cerutiUnconventionalRobustIntegrator2021}.

In the numerical analysis of kinetic equations, the dynamical low-rank method has recently 
been applied to many problems. Here we mention a few representative ones:
the Vlasov equation
\cite{einkemmerLowRankProjectorSplittingIntegrator2018},
\cite{einkemmerMassMomentumEnergy2021a},
linear transport equation
\cite{einkemmerAsymptoticpreservingDynamicalLowrank2021},
\cite{pengLowrankMethodTwodimensional2020},
\cite{dingDynamicalLowRankIntegrator2021},
Boltzmann equation
\cite{huAdaptiveDynamicalLow2021a},
and BGK equation
\cite{einkemmerLowRankAlgorithmWeakly2019}
\cite{einkemmerEfficientDynamicalLowRank2021}.
In particular, the last contribution \cite{einkemmerEfficientDynamicalLowRank2021} is significant for preserving the asymptotic fluid limit of the
collisional BGK equation, which inspires our current work.

In this paper we present a dynamical low-rank algorithm for the solution of a
model equation for the Vlasov equation with collisions. Collisions with a
Coulomb interaction potential can be described by an integro-differential
operator with a drift term and a diffusion term, i.e., the Landau operator
\cite{goldstonIntroductionPlasmaPhysics1995},
\cite{cercignaniBoltzmannEquationIts1988},
\cite{rosenbluthFokkerPlanckEquationInverseSquare1957},
\cite{villaniReviewMathematicalTopics2002}
or simplified Fokker-Planck type operator
\cite{doughertyModelFokkerPlanckEquation1964}.
We therefore consider
the non-magnetic Vlasov equation with a linear Fokker-Planck collision
operator, in the so-called ``high-field'' limit. This scaling was introduced in
\cite{poupaudRunawayPhenomenaFluid1992} as a model for the semiconductor
Boltzmann equation, was treated numerically in
\cite{cercignaniDeviceBenchmarkComparisons}, and in
\cite{jinAsymptoticPreservingScheme2011} with an asymptotic-preserving scheme.
It retains key properties of the full Vlasov-Landau-Poisson equation,
specifically the diffusive collision operator and nonlinear coupling between $f$
and the electric field. As such, it provides an interesting test case for the
dynamical low-rank method applied to collisional plasma equations. In a
nondimensionalized form, the equation reads
\begin{align}
\label{eqn:vpfp_original}
\partial_t f + v \cdot \nabla_x f + \frac{1}{\epsilon} E \cdot \nabla_v f = \frac{1}{\epsilon} P_{\text{FP}}(f), \qquad t > 0, \quad x \in \Omega \subset R^{d}, \quad v \in \mathbb{R}^{d}.
\end{align}
The function $f(x, v,t)$ is the single-particle probability density function
defined over phase space of $d$ spatial and $d$ velocity dimensions. The
operator $P_{\text{FP}}$ is the linear Fokker-Planck operator
\begin{equation}
P_{\text{FP}}(f) = \nabla_v \cdot (vf + \nabla_v f).
\end{equation}
One can imagine that $f$ describes a population of electrons moving under the influence
of their own inertia and electric field. The small parameter $\epsilon > 0$ is a
scaling parameter which determines the strength of both collisions and the
electric field $E$. It is analogous to the Knudsen number from the theory of
hydrodynamic limits of the Boltzmann equation.

The electric field $E$ is determined self-consistently from the phase space
distribution $f$ via Amp\`{e}re's law:
\begin{equation}
\label{eqn:ampere}
\partial_{t} E = -J, \qquad t > 0, \quad x \in \Omega,
\end{equation}
where the current density $J$ is defined as
\begin{equation}
  J(x,t) = \langle v,f \rangle_{v}.
\end{equation}
Here we have made use of a notation for the $L^{2}$ inner product, which we now
define as
\begin{equation}
  \langle g, h \rangle_{x} = \int_{\Omega} gh \rd{x}, \quad \langle g, h \rangle_{v} = \int_{\mathbb{R}^{d}} gh \rd{v}, \quad \langle g, h \rangle_{xv} = \int_{\Omega} \int_{\mathbb{R}^{d}} gh \rd{x}\rd{v}.
\end{equation}
The initial electric field will be specified via a static background charge
density $\eta(x)$. To continue the physical picture of an electron fluid, $\eta$
may represent a density of ions which do not move on the timescale resolved by
(\ref{eqn:vpfp_original}). To be physical, the field $E$ should satisfy Gauss's law with respect to the density $\rho$:
\begin{equation}
    \label{eqn:gauss}
    E = -\nabla_{x} \phi(x), \quad -\nabla_{x}^{2} \phi(x) = \rho(x,t) - \eta(x), \\
\end{equation}
where
\begin{equation}
\rho(x,t) = \langle 1,f \rangle_{v}.
\end{equation}
It is easy to show that if $E(x, t)$ satisfies Amp\`{e}re's equation (\ref{eqn:ampere}) and
satisfies Gauss's law (\ref{eqn:gauss}) at time 0, then it will satisfy
(\ref{eqn:gauss}) for all time. Numerically, we initialize $E(x, 0)$ using
Gauss's law and a specified background density $\eta(x)$, and then timestep $E$
using Amp\`{e}re's law. This introduces an error in Gauss's law which grows in time.
Codes which care about the detailed electrostatic properties of such systems
must use ``divergence-cleaning'' methods to clear this error; we will simply
note that it exists but is first-order in time.

To recapitulate, in this work we are solving the coupled system
\begin{align}
  \label{eqn:full_system}
\begin{cases}
\partial_tf + v\cdot \nabla_x f + \frac{1}{\epsilon} E \cdot \nabla_{v} f = \frac{1}{\epsilon} \nabla_v \cdot (vf + \nabla_{v} f), & (x, v, t) \in \Omega \times \mathbb{R}^{d} \times [0, T] \\
\partial_tE = -J, & (x, t) \in \Omega \times [0, T] \\
E = -\nabla_{x} \phi, \quad  -\nabla_{x}^{2} \phi = \rho - \eta, & x \in \Omega,\ t = 0.
\end{cases}
\end{align}

\subsection{Asymptotic fluid limit}
The limit of (\ref{eqn:full_system}) for very small $\epsilon$ is a sort of
electrostatic ``creeping flow'', in which inertial forces are vanishingly small
compared to electrostatic forces.
To analyze the limit \(\epsilon \rightarrow 0\), we introduce a scaled ``local
Maxwellian'' defined as
\begin{equation} \label{Maxwellian}
M (x,v,t)= \frac{1}{(2\pi)^{d/2}} e^{-\frac{|v - E(x,t)|^2}{2}}.
\end{equation}
One should note that, in contrast to the Maxwellian equilibrium of the Boltzmann
equation, this function has a uniform density. It is isothermal and the flow
velocity is equal to $E$.
It is not hard to show that (\ref{eqn:vpfp_original}) is equivalent to
\begin{align}
\label{eqn:vpfp_combined}
\partial_tf + v\cdot \nabla_x f = \frac{1}{\epsilon} \nabla_v \cdot \left[ M \nabla_v (M^{-1} f) \right].
\end{align}
To see this we expand the right hand side of (\ref{eqn:vpfp_combined}):
\begin{align}
\nabla_v (M^{-1}f) &= (\nabla_v M^{-1})f + M^{-1}\nabla_v f = (v - E)M^{-1}f + M^{-1}\nabla_v f= M^{-1} (v - E + \nabla_v)f.
\end{align}
Therefore,
\begin{align}
\nabla_v \cdot [M \nabla_v (M^{-1}f)] &= \nabla_v \cdot \left[ (v - E + \nabla_v) f \right]= \underbrace{\nabla_v \cdot (vf + \nabla_v f)}_{P_{\text{FP}}(f)} - \cancel{(\nabla_v \cdot E) f} - \underbrace{E \cdot \nabla_v f}_{\text{force term}},
\end{align}
where we have accounted for both \(P_{\text{FP}}(f)\) and the force term on the left
hand side of (\ref{eqn:vpfp_original}).

The form (\ref{eqn:vpfp_combined}) reveals the dominant balance structure of
(\ref{eqn:vpfp_original}): the linear Fokker-Planck operator and the
electrostatic force term are balanced under this scaling; together they derive $f$ to a local equilibrium. Indeed, when $\epsilon \rightarrow 0$ in (\ref{eqn:vpfp_combined}), formally $\nabla_v \cdot \left[ M \nabla_v (M^{-1} f) \right] \rightarrow 0$ which implies $f\rightarrow \rho M$. To see this, first note that
\begin{equation}
\begin{split}
\int \nabla_v \cdot  [M \nabla_v (M^{-1}f)] \log (M^{-1} f)\rd{v}&=\int \nabla_v \cdot  [f  \nabla_v \log (M^{-1}f)]  \log (M^{-1} f)\rd{v}\\
&=-\int  f \left| \nabla_v \log (M^{-1}f)\right|^2\rd{v}\leq 0.
\end{split}
\end{equation}
Accordingly one can derive (by a cycle of implications)
\begin{equation}
\int \nabla_v \cdot  [M \nabla_v (M^{-1}f)] \log (M^{-1} f)\rd{v}=0 \Longleftrightarrow f=cM \Longleftrightarrow \nabla_v \cdot  [M \nabla_v (M^{-1}f)] =0,
\end{equation}
where $c$ is a function independent of $v$. Finally $\rho=\langle 1,f\rangle_v=c\langle1, M\rangle_v=c$.

To derive a macroscopic system of (\ref{eqn:vpfp_combined}) when $\epsilon \rightarrow 0$, we first take moments $\langle 1,\cdot \rangle_v$, $\langle v,\cdot \rangle_v$ of (\ref{eqn:vpfp_combined}) to obtain
\begin{align}
&\partial_t\rho+\nabla_x\cdot J=0, \label{moment-1}\\
&\partial_t J+\nabla_x\cdot \langle v\otimes v,f \rangle_v=\frac{1}{\epsilon}(\rho E-J).\label{moment-2}
\end{align}
As $\epsilon\rightarrow 0$, one has $J\rightarrow \rho E$ from (\ref{moment-2}).
Then (\ref{moment-1}) becomes
\begin{equation}
  \label{eqn:mass_conservation}
\partial_t\rho+\nabla_x\cdot (\rho E)=0,
\end{equation}
which together with Amp\`{e}re's law (\ref{eqn:ampere}) constitute the limiting system:
\begin{equation}
\label{eqn:limiting_system}
\begin{cases}
  \partial_{t} \rho + \nabla_{x} \cdot (\rho E) = 0, & (x, t) \in \Omega \times [0, T] \\
  \partial_tE = -\rho E, & (x, t) \in \Omega \times [0, T] \\
 E = -\nabla_{x} \phi, \quad  -\nabla_{x}^{2} \phi = \rho - \eta, & x \in \Omega,\ t = 0.
\end{cases}
\end{equation}

The system (\ref{eqn:limiting_system}) fully determines the behavior of the
kinetic system (\ref{eqn:full_system}) in the asymptotic limit
$\epsilon \rightarrow 0$. Our numerical scheme is careful to preserve this
asymptotic limit at the discrete level. However, trying to design a low-rank
scheme that smoothly approaches this limit quickly runs into a problem:
$f^{0} = \rho M$ is not a low-rank function in $x$ and $v$, i.e., we cannot
write it in the form (\ref{eqn:lowrank_ansatz}) with small $r$. Recall the
definition (\ref{Maxwellian}), where the cross term $e^{-v \cdot E(x, t)}$ is not
low rank. If we require a high rank to resolve the limiting solution $f^{0}$,
then we are, in a sense, wasting effort on a kinetic system whose dominant dynamics are
described by the much lower-dimensional system (\ref{eqn:limiting_system}).

To resolve this problem, we can observe that while $f^{0} = \rho M$ is not low
rank in $x$ and $v$, $\rho(x,t)$ certainly is. This motivates us to consider a
low-rank approximation to the quotient
\begin{equation}
  g(x, v, t) = M^{-1}f(x, v, t),
\end{equation}
which as we have seen has a rank-1 asymptotic limit.
We will therefore search for solutions of the form
\begin{equation}
  f = Mg = \frac{1}{(2\pi)^{d/2}} e^{-\frac{|v-E(x,t)|^{2}}{2}} g,
\end{equation}
where \(g\) is given a low-rank approximation $\tilde{g}$:
\begin{equation}
\label{eqn:g_low_rank}
g(x, v, t) \approx \tilde{g}(x, v, t) := \sum_{i, j = 1}^r X_i(x, t) S_{ij}(t) V_j(v, t).
\end{equation}
The bases $X_{i}$ and $V_{j}$ are required to satisfy orthogonality relations,
\begin{equation}
  \label{eqn:orthogonality}
  \langle X_{i}, X_{k} \rangle_{x} = \delta_{ik}, \quad \langle V_{j}, V_{l} \rangle_{v} = \delta_{jl}.
\end{equation}
The approximation $\tilde{g}$ is the quantity which we will timestep using the
dynamical low-rank method. We will also solve Amp\`{e}re's law to advance the
electric field, which has only a dependence on $x$. From these, we can
reconstruct the approximate solution $\tilde{f} = \tilde{g}M$ as desired.

The rest of this paper is organized as follows. In Section 2 we present the
dynamical low-rank algorithm for the evolution of the function $\tilde{g}$
defined in (\ref{eqn:g_low_rank}). This consists of deriving the PDEs satisfied by
the low-rank factors. In Section 3 we present a first-order time integration
scheme for that system of PDEs. In Section 4, we address the question of
discretization in physical ($x$) and velocity ($v$) space. Section 5 consists of
a brief discussion of the asymptotic limit of the discrete system derived in
Sections 2-4, verifying that it recovers the fluid equations of
(\ref{eqn:limiting_system}). Finally, Section 6 includes some numerical results
from an implementation of our algorithm and comparison with the full tensor
solution of the kinetic equation.

\section{Dynamical low-rank algorithm}

The dynamical low-rank algorithm works by confining the time derivative of the
system to the tangent space of a low-rank manifold. We will write down the
time derivative of $g$ imposed by (\ref{eqn:vpfp_combined}), and then discuss
its projection onto the tangent space.
We can derive the dynamics for $g$ by plugging \(f = Mg\) into
(\ref{eqn:vpfp_combined}). This gives
\begin{align}
\partial_t (Mg) + v\cdot \nabla_x (Mg) &= \frac{1}{\epsilon} \nabla_v \cdot [M \nabla_v g] \nonumber \\
\Longrightarrow \ g \partial_t M + M \partial_t g + v\cdot (g \nabla_x M + M \nabla_x g) &= \frac{1}{\epsilon} (\nabla_v M \cdot \nabla_v g) + \frac{1}{\epsilon} M\nabla_v^2 g.
\end{align}
Consolidating terms and dividing through by $M$, we obtain
\begin{align}
\label{eqn:h_defn}
\partial_t g &= -v \cdot \nabla_x g + \frac{1}{\epsilon} \left[ (\nabla_v - v + E) \cdot \nabla_v g \right] - \frac{1}{M} (\partial_t M + v\cdot \nabla_x M) g \nonumber \\
&= -v \cdot \nabla_x g + \frac{1}{\epsilon} \left[ (\nabla_v - v + E) \cdot \nabla_v g \right] - \mathcal{M} g \\
&:= \mathcal{H}[g], \nonumber
\end{align}
where we have introduced the shorthand
\begin{equation}
  \label{eqn:mathcalM}
  \mathcal{M} = \frac{1}{M}(\partial_{t}M + v\cdot \nabla_{x} M).
\end{equation}
The time derivative of the low-rank approximation is now given by composing
$\mathcal{H}$ with a projection operator. That is,
\begin{equation}
  \partial_{t} \tilde{g} = P(\mathcal{H}[\tilde{g}]),
\end{equation}
where $P$ is the projection onto the tangent space to the manifold of functions
with a rank $r$ representation as in (\ref{eqn:g_low_rank}). It can be shown
(\cite{kochDynamicalLowRank2007},
\cite{lubichProjectorsplittingIntegratorDynamical2014}) that the projection
operator takes the form
\begin{equation}
  P(h) = \sum_{j} \langle V_{j}, h \rangle_{v} V_{j} - \sum_{ij} X_{i} \langle X_{i} V_{j}, h \rangle_{xv} + \sum_{i} X_{i} \langle X_{i}, h \rangle_{x}.
\end{equation}
We now have a three-term expression for the time derivative of our
low-rank approximation:
\begin{equation}
  \label{eqn:dt_g_unsplit}
  \partial_{t} \tilde{g} = \sum_{j} \langle V_{j}, \mathcal{H}[\tilde{g}] \rangle_{v} V_{j} - \sum_{ij} X_{i} \langle X_{i} V_{j}, \mathcal{H}[\tilde{g}] \rangle_{xv}V_{j} + \sum_{i} X_{i} \langle X_{i}, \mathcal{H}[\tilde{g}] \rangle_{x}.
\end{equation}
This form lends itself to a first-order-in-time Lie-Trotter operator splitting,
which we will employ in this paper. Higher-order splitting schemes are possible,
for example a second-order scheme based on Strang splitting, although this
requires extra care to properly center the electric field
\cite{einkemmerLowRankProjectorSplittingIntegrator2018}. The first-order-in-time
scheme splits (\ref{eqn:dt_g_unsplit}) into the three equations
\begin{align}
  \label{eqn:g_split_1}
  \partial_{t} \tilde{g} &= \sum_{j} \langle V_{j}, \mathcal{H}[\tilde{g}] \rangle_{v} V_{j}, \\
  \label{eqn:g_split_2}
  \partial_{t} \tilde{g} &= - \sum_{ij} X_{i} \langle X_{i} V_{j}, \mathcal{H}[\tilde{g}] \rangle_{xv} V_{j}, \\
  \label{eqn:g_split_3}
  \partial_{t} \tilde{g} &= \sum_{i} X_{i} \langle X_{i}, \mathcal{H}[\tilde{g}] \rangle_{x}.
\end{align}
We actually implement this scheme in terms of a pair auxiliary bases (making use
of (\ref{eqn:orthogonality})),
\begin{align}
  K_{j}(x, t) &= \langle \tilde{g}, V_{j} \rangle_{v} = \sum_{i} X_{i}(x, t)S_{ij}(t),\\
L_{i}(v, t) &= \langle X_{i}, \tilde{g} \rangle_{x} = \sum_{j} S_{ij}(t) V_{j}(v, t).
\end{align}
With this notation the splitting scheme is as follows:
\begin{itemize}
\item The first step holds the $V_{j}$ basis constant.
Take the inner product of (\ref{eqn:g_split_1}) with $V_{j}$ to obtain
\begin{equation}
  \label{eqn:K_split}
  \partial_{t} K_{j} = \langle V_{j}, \mathcal{H}[\tilde{g}] \rangle_{v}.
\end{equation}
Integrate this equation for one time step to obtain a new value for $K_{j}$.
Then perform a QR decomposition of $K_{j}$ to obtain a new orthogonal basis
$X_{i}$ and coefficients $S_{ij}$.
\item The second step holds both bases constant. Take the inner
product of (\ref{eqn:g_split_2}) with $V_{j}$ in $v$, and with $X_{i}$ in $x$ to obtain
\begin{equation}
  \label{eqn:S_split}
  \partial_{t} S_{ij} = -\langle X_{i}V_{j}, \mathcal{H}[\tilde{g}] \rangle_{xv}.
\end{equation}
Integrate this equation for one time step to obtain a new matrix $S_{ij}$.
\item The third step holds the $X_{i}$ basis constant. Take the inner product of
(\ref{eqn:g_split_3}) with $X_{i}$ to obtain
\begin{equation}
  \label{eqn:L_split}
  \partial_{t} L_{i} = \langle X_{i}, \mathcal{H}[\tilde{g}] \rangle_{x}.
\end{equation}
Integrate this equation for one time step to obtain a new value for $L_{i}$.
Then perform a QR decomposition of $L_{i}$ to obtain a new orthogonal basis
$V_{j}$ and coefficient matrix $S_{ij}$.
\end{itemize}
The above algorithm has the excellent property that it is robust to
``overapproximation'', i.e. small singular
values in \(S\) \cite{lubichProjectorsplittingIntegratorDynamical2014}.

The time splitting scheme for $\tilde{g}$ may be straightforwardly coupled with Amp\`{e}re's equation
(\ref{eqn:ampere}), which can be written in terms of the low-rank components as
\begin{equation}
  \label{eqn:ampere_lowrank}
  \partial_{t} E = -\langle v M \tilde{g} \rangle_{v} = -\sum_{i,j} X_{i} S_{ij} \langle v M V_{j} \rangle_{v}.
\end{equation}

\subsection{Time evolution of low-rank components}

In this section we expand the inner products involving $\mathcal{H}[\tilde{g}]$ which
appear in equations (\ref{eqn:K_split}), (\ref{eqn:S_split}),
(\ref{eqn:L_split}). The result is a self-contained system of $r$ coupled PDEs
for $K_{j}$ and $L_{i}$, and a matrix-valued ODE for $S_{ij}$ of size $r \times r$.

Plugging (\ref{eqn:h_defn}) and (\ref{eqn:g_low_rank}) into (\ref{eqn:K_split}) gives
\begin{align}
\partial_t K_j &= -\sum_{k, l} \langle V_j, v\cdot(\nabla_x X_k) S_{kl} V_l \rangle_v - \sum_{kl} X_k S_{kl} \langle V_j, V_l \mathcal{M} \rangle_v \nonumber \\
&\quad\quad + \frac{1}{\epsilon} \left( \sum_{kl} X_k S_{kl} \langle V_j [(\nabla_v - v + E) \cdot \nabla_v V_l] \rangle_v \right) \nonumber\\
\label{eqn:partial_t_K_expanded}
&= -\sum_{l} (\nabla_x K_l) \cdot \langle v V_j V_l \rangle_v - \sum_{l} K_l \langle V_j V_l \mathcal{M} \rangle_v \\
&\quad\quad + \frac{1}{\epsilon} \left( \sum_{l} K_l \left[ \langle V_j (\nabla_v - v)\cdot \nabla_v V_l \rangle_v + E \cdot \langle V_j \nabla_v V_l \rangle_v \right] \right). \nonumber
\end{align}

Plugging (\ref{eqn:h_defn}) and (\ref{eqn:g_low_rank}) into (\ref{eqn:S_split}) gives
\begin{align}
\partial_t S_{ij} &= \sum_{kl} \langle X_i S_{kl} \nabla_x X_k \cdot \langle v V_j V_l \rangle_v \rangle_x + \sum_{kl} \langle X_i X_k S_{kl} V_l V_j \mathcal{M} \rangle_{xv} \nonumber\\
&\quad\quad -\frac{1}{\epsilon} \left( \sum_{kl} \left\langle X_i S_{kl} X_k \langle V_j [\nabla_v - v + E] \cdot \nabla_v V_l \rangle_v \right\rangle_x \right) \nonumber\\
\label{eqn:partial_t_S_expanded}
&= \sum_{kl} S_{kl} \langle X_i \nabla_x X_k \rangle_x \cdot \langle v V_j V_l \rangle_v + \sum_{kl} S_{kl} \langle X_i X_k V_l V_j \mathcal{M} \rangle_{xv} \\
&\quad\quad -\frac{1}{\epsilon} \left( \sum_{kl} S_{kl} \left[\langle X_i X_k \rangle_x \langle V_j (\nabla_v - v) \cdot \nabla_v V_l \rangle_v +  \langle X_i X_k E \rangle_x \cdot \langle V_j \nabla_v V_l \rangle_v\right] \right). \nonumber
\end{align}

Plugging (\ref{eqn:h_defn}) and (\ref{eqn:g_low_rank}) into (\ref{eqn:L_split}) gives
\begin{align}
\partial_t L_i &= -\sum_{kl} v\cdot \langle X_i (\nabla_x X_k) S_{kl} V_l \rangle_x - \sum_{kl} \langle X_i X_k S_{kl} V_l \mathcal{M} \rangle_x \nonumber\\
&\quad\quad+ \frac{1}{\epsilon} \left( \sum_{kl} \langle X_i X_k S_{kl} [(\nabla_v - v + E) \cdot \nabla_v V_l] \rangle_x \right) \nonumber\\
\label{eqn:partial_t_L_expanded}
&= -\sum_{k} v\cdot \langle X_i (\nabla_x X_k) \rangle_x L_k - \sum_{k} \langle X_i X_k  \mathcal{M} \rangle_x L_k \\
&\quad\quad+ \frac{1}{\epsilon} \left(\sum_{k} \left[ \langle X_i X_k \rangle_x (\nabla_v - v) + \langle X_i X_k E \rangle_x \right] \cdot \nabla_v L_k\right). \nonumber
\end{align}

We also expand the terms involving \(\mathcal{M}\) (defined in (\ref{eqn:mathcalM})):
\begin{equation*}
  \partial_{t} M = (v - E) \cdot (\partial_{t} E) M= - M (v - E) \cdot J,
\end{equation*}
\begin{equation*}
  v \cdot \nabla_{x} M = M \sum_{i,j} (v_{j}-E_{j})v_{i}\partial_{x_{i}} E_{j} = M \sum_{i,j} (v_{j}v_{i}\partial_{x_{i}} E_{j} - E_{j} v_{i} \partial_{x_{i}} E_{j}) = M(v \otimes v):\nabla_{x} E - \frac{M}{2} v\cdot \nabla_{x}(E^{2}).
\end{equation*}
Note that $\frac{1}{2}\nabla_x (E^2)\neq E\cdot \nabla_x E$. $(\nabla_xE)_{ij}:=\partial_{x_j} E_i$ and $A:B:=\sum_{ij}a_{ij}b_{ij}.$
Putting these together we obtain
\begin{align}
\label{eqn:mathcal_M}
\mathcal{M} &= \frac{1}{M}(\partial_t M + v\cdot\nabla_x M) = E \cdot J - v\cdot J - \frac{1}{2}v\cdot\nabla_{x}(E^{2}) + (v\otimes v):\nabla_{x} E \nonumber \\
&:=  \Mone + v\cdot\Mtwo + (v \otimes v):\Mthree,
\end{align}
where
\begin{equation}
\label{eqn:numbered_Mcal}
\Mone = E\cdot J, \quad \Mtwo = -J - \frac{1}{2}\nabla_x (E^2), \quad \Mthree = \nabla_x E.
\end{equation}
Here we use $\textbf{boldface}$ to denote vectors of length $d$, and
$\textsf{sans-serif}$ to denote tensors of size $d \times d$. Both vectors and
tensors may also vary in $x$ and $v$. In all cases the tensor
contractions $\cdot, :$ indicate contraction over the length-$d$ dimensions. The
terms involving \(\mathcal{M}\) then expand to
\begin{align*}
\langle V_j V_l \mathcal{M} \rangle_v &= \delta_{jl} \Mone + \langle v V_j V_l \rangle_v \cdot \Mtwo + \langle (v \otimes v) V_j V_l \rangle_v : \Mthree, \\
\langle X_i X_k \mathcal{M} \rangle_x &= \langle X_i X_k \Mone \rangle_x + v \cdot \langle X_i X_k \Mtwo \rangle_x + (v \otimes v) : \langle X_i X_k \Mthree \rangle_x, \\
\langle X_i X_k V_j V_l \mathcal{M} \rangle_{xv} &= \delta_{jl} \langle X_i X_k \Mone \rangle_x + \langle v V_j V_l \rangle_v \cdot \langle X_i X_k \Mtwo \rangle_x \\
&\quad\quad+ \langle (v \otimes v) V_j V_l \rangle_v : \langle X_i X_k \Mthree \rangle_x.
\end{align*}

\section{First order in time scheme}

The algorithm described up to this point has been fully continuous, except for
the projection onto the low-rank manifold. We now present a discretization in
time, leaving space continuous for the moment. The time discretization makes use
of an implicit-explicit (IMEX) scheme for capturing the fast dynamics of the
collision operator in the fluid limit ($\epsilon \ll 1$).

In the following we report rough estimates of the computational complexity of each substep.
To avoid complicating the presentation unnecessarily, for these estimates we consider $d \sim 1$,
so that we are free to ignore both the dimension and constant factors in our ``big-O'' notation.

Suppose we have the quantities \((E^n, X_i^n, V_j^n, S_{ij}^n)\) at timestep
\(t^n\). Then we calculate \((E^{n+1},X_i^{n+1},
V_j^{n+1}, S_{ij}^{n+1})\) in the following way.

\subsection{Step 1: Update $E$}
\begin{enumerate}
\item Compute the following integral appearing in (\ref{eqn:ampere_lowrank}):
\begin{align}
I_j^n(x) := \langle v V_j^n M^n \rangle_v = \frac{1}{(2\pi)^{d/2}}\int v V_j^n(v) e^{-\frac{|v - E^n(x)|^2}{2}}\rd{v}.
\end{align}
A naive computation of this integral requires $\bO(N_{x}N_{v})$ steps, a
computational cost that is unacceptably high. However, since the Maxwellian is
isothermal, the integral has a convolutional structure, and may be computed
using a Fast Fourier Transform (FFT). The required substeps are:
\begin{itemize}
\item Compute the convolution
\begin{align}
\label{eqn:Ijn_conv}
\ell^n_j(\zeta) = \left[ (v \mapsto v V_j^n(v)) * \left( v \mapsto e^{-|v|^2/2} \right) \right](\zeta)
\end{align}
using an FFT.

\textbf{Cost:} $\bO(rN_{v}\log N_{v})$.
\item Compute the composition of \(E^n\) with \(\ell^n_j\) using any interpolation
scheme from the FFT nodes to an arbitrary point \(E^{n}(x)\):
\begin{align*}
I_j^n(x) = \frac{1}{(2\pi)^{d/2}}\ell^n_j(E^n(x)).
\end{align*}
\textbf{Cost:} $\bO(N_{x})$.
\end{itemize}
Exploiting the convolutional structure of $I_{j}^{n}(x)$ with an FFT reduces the total
computational cost to $\bO(rN_{v} \log{N_{v}})$, which is acceptable.
\item Compute the current density:
\begin{align}
J^n(x) = \sum_{ij} X_i^n(x) S_{ij}^n I_j^n(x).
\end{align}
\textbf{Cost:} $\bO(r^{2}N_{x})$.

\item Perform a Forward Euler step to solve (\ref{eqn:ampere_lowrank}):
\begin{align}
  E^{n+1}(x) &= E^n(x) - \Delta t J^n.
\end{align}
\textbf{Cost:} $\bO(N_{x})$.
\end{enumerate}

\subsection{Step 2: Update \(X\), \(S\), and \(V\)}
\subsubsection{$K$ step}
\begin{enumerate}
  \item Compute integrals in \(v\). We use $\textbf{boldface}$ to denote vector-valued matrices of total size $r \times r \times d$,
  and $\textsf{sans-serif}$ to denote tensor-valued matrices of total size $r \times r \times d \times d$. In both cases the indices running over the length-$d$ dimensions are suppressed.
  The integrals to compute are:
\begin{align}
\label{eqn:defn_c_ints}
\mathbf{c}^1_{jl} = \langle v V^n_j V^n_l \rangle_v, \quad \overline{\overline{\mathsf{c}}}^2_{jl} = \langle (v \otimes v) V^n_j V^n_l \rangle_v,
\end{align}
\begin{align}
\label{eqn:defn_d_ints}
d^1_{jl} = \langle V^n_j (\nabla_v - v)\cdot \nabla_v V^n_l \rangle_v, \quad \mathbf{d}^2_{jl} = \langle V^n_j \nabla_v V^n_l \rangle_v.
\end{align}
\textbf{Cost:} $\bO(r^{2}N_{v})$.
\item Compute $\Mone^{n}, \Mtwo^{n}$, and $\Mthree^{n}$:
\begin{align}
  \Mone^{n} = E^n \cdot J^n, &\quad \Mtwo^{n} = -J^n - \frac{1}{2}\nabla_x((E^n)^2), \\
  &\Mthree^{n} = \nabla_x E^n.
\end{align}
\textbf{Cost:} $\bO(N_{x})$.
\item Compute matrices on the right hand side of (\ref{eqn:partial_t_K_expanded}).
\begin{align}
\label{eqn:K_matrix_A}
A^1_{jl} &= \delta_{jl} \Mone^{n} + \mathbf{c}_{jl}^1 \cdot \Mtwo^{n} + \overline{\overline{\mathsf{c}}}^2_{jl} : \Mthree^{n}, \\
A^2_{jl} &= d^1_{jl} + E^{n} \cdot \mathbf{d}^2_{jl}.
\end{align}
\textbf{Cost:} $\bO(r^{2}N_{x})$.

\item The evolution equation (\ref{eqn:partial_t_K_expanded}) for \(K\) may now be written as
\begin{align}
\label{eqn:partial_t_K_matrix_form}
\partial_t K_j &= -\sum_{l} \mathbf{c}_{jl}^1 \cdot \nabla_x K_l - \sum_{l} A_{jl}^1 K_l + \frac{1}{\epsilon} \sum_{l} A_{jl}^2 K_l.
\end{align}
Advance (\ref{eqn:partial_t_K_matrix_form}) in time, using an IMEX step to handle the stiff
term:
\begin{align}
\label{eqn:K_advance}
\sum_{l} \left[ \delta_{jl} - \frac{\Delta t}{\epsilon} A_{jl}^2 \right] K_l^{n+1} = K_j^n - \Delta t \left( \sum_{l} \mathbf{c}_{jl}^1 \cdot \nabla_x K_l^n + \sum_{l} A_{jl}^1 K_l^n \right).
\end{align}
Note that the only differential operator, namely $\nabla_{x}$, appearing
in this equation is treated explicitly. Therefore the linear system
appearing in this equation involves no coupling between points in $x$.
When discretized it will consist of $N_{x}$ separate systems each of
size $r \times r$. We can solve this small system at each point in $x$
using any standard dense linear solver---the size is not large enough to
warrant any special technique.

\textbf{Cost:} $\bO(r^{2}N_{x})$ for both the right-hand side and the implicit step, due to solving each $r \times r$ system separately.

\item Perform a QR decomposition of \(K^{n+1}_j\) to obtain \(X^{n+1}_i\) and \(S^1_{ij}\).

\textbf{Cost:} $\bO(r^{2}N_{x})$.  
\end{enumerate}

\subsubsection{$S$ step}
\label{sec:org27f4e10}
\begin{enumerate}
\item Compute the integrals in $x$, using the new basis $X^{n+1}$:
\begin{align}
\label{eqn:defn_c_ints_X}
&c^\star_{ik} = \langle X^{n+1}_i X^{n+1}_k \Mone^{n} \rangle_x, \quad \mathbf{c}^{\star\star}_{ik} = \langle X^{n+1}_i X^{n+1}_k \Mtwo^{n} \rangle_x, \quad \overline{\overline{\mathsf{c}}}^{\star\star\star}_{ik} = \langle X^{n+1}_i X^{n+1}_k \Mthree^{n} \rangle_x, \\
&\mathbf{d}^\star_{ik} = \langle X^{n+1}_i \nabla_x X^{n+1}_k \rangle_x, \quad \mathbf{e}^\star_{ik} = \langle X^{n+1}_i X^{n+1}_k E^n \rangle_x.
\end{align}
\textbf{Cost:} $\bO(r^{2} N_{x})$.
\item Compute
\begin{align}
\tilde{c}_{ij;kl} = \langle X^{n+1}_i X^{n+1}_k V^n_j V^n_l \mathcal{M}^n \rangle_{xv} = \delta_{jl} c^\star_{ik} + \mathbf{c}^1_{jl} \cdot \mathbf{c}^{\star\star}_{ik} + \overline{\overline{\mathsf{c}}}^2_{jl} : \overline{\overline{\mathsf{c}}}^{\star\star\star}_{ik}.
\end{align}
\textbf{Cost:} $\bO(r^{4})$.
\item Compute the order-four tensors
\begin{align}
\label{eqn:defn_B}
B^{1}_{ij;kl} &= d^{\star}_{ik}\cdot\mathbf{c}^{1}_{jl} + \tilde{c}_{ij;kl}, \\
B^{2}_{ij;kl} &= \delta_{ik}d^1_{jl} + \mathbf{e}^\star_{ik} \cdot \mathbf{d}^2_{jl}.
\end{align}
\textbf{Cost:} $\bO(r^{4})$.

\item The evolution equation for $S$ may now be written as
\begin{equation}
  \label{eqn:partial_t_S}
  \partial_{t} S_{ij} = \sum_{kl} B^{1}_{ij;kl} S_{kl} - \frac{1}{\epsilon} \sum_{kl} B^{2}_{ij;kl} S_{kl}.
\end{equation}
Perform a Forward Euler step to advance \(S^1_{ij} \rightarrow S^2_{ij}\):
\begin{align}
\label{eqn:S_forward_euler}
S^2_{ij} = S^1_{ij} + \Delta t \sum_{kl} S^1_{kl} B^{1}_{ij;kl} - \frac{\Delta t}{\epsilon} \sum_{kl} S^1_{kl} B^{2}_{ij;kl}.
\end{align}
Our use of a Forward Euler step here differs from the presentation in
\cite{einkemmerEfficientDynamicalLowRank2021}, where an IMEX step was
used to advance the $S$ equation in the case of the BGK collision operator. Since the structure of the Fokker-Planck operator is more complicated than the BGK type, special care is needed. Heuristically, one can see that something different from the $K$ and $L$ steps may be required, simply because the S equation runs backwards in time. For a more detailed
justification and a discussion of how the situation differs from \cite{einkemmerEfficientDynamicalLowRank2021}, refer to Appendix \ref{S_STEP}.

\textbf{Cost:} $\bO(r^{4})$.
\end{enumerate}

\subsubsection{$L$ step}
\label{sec:orgbfbe5fc}
\begin{enumerate}
\item Compute the $r \times r$ matrix
\begin{align}
\hat{c}_{ik} = c^\star_{ik} + v\cdot \mathbf{c}^{\star\star}_{ik} + (v \otimes v) : \overline{\overline{\mathsf{c}}}^{\star\star\star}_{ik}.
\end{align}
\textbf{Cost:} $\bO(r^{2})$.

\item The equation (\ref{eqn:partial_t_L_expanded}) may now be written as
\begin{align}
  \label{eqn:partial_t_L_discrete}
\partial_t L_i = -\sum_{k} v \cdot \mathbf{d}^\star_{ik} L_k - \sum_{k}\hat{c}_{ik} L_k + \frac{1}{\epsilon} \sum_{k}\left( \delta_{ik} \nabla_v^2 - \delta_{ik} v \cdot \nabla_{v} + \mathbf{e}^\star_{ik} \cdot \nabla_v \right) L_k.
\end{align}
Advance (\ref{eqn:partial_t_L_discrete}) using an IMEX step by solving the system
\begin{align}
\label{eqn:L_IMEX}
\sum_{k} \left[ \delta_{ik} - \frac{\Delta t}{\epsilon} (\delta_{ik} \nabla^2_v - \delta_{ik} v \cdot \nabla_{v} + \mathbf{e}^\star_{ik}\cdot \nabla_v ) \right] L^{n+1}_k = L_i^n - \Delta t \sum_{k} (v \cdot \mathbf{d}^\star_{ik} + \hat{c}_{ik}) L_k^n.
\end{align}
Note that in contrast to (\ref{eqn:K_advance}), the left-hand side of
this equation does involve differential operators in $v$, and so the
linear system may be discretized by a fully coupled (but sparse) matrix
of size $rN_{v} \times rN_{v}$. There are $\bO(rN_{v})$ non-empty
entries. Using Krylov subspace methods lets us keep the total cost of solving this system
on the order of $\bO(rN_{v})$, assuming the number of iterations does not
grow unboundedly with $r$ or $N_{v}$, which is what we observe in practice.

\textbf{Cost:} $\bO(r^{2}N_{v})$.

\item Perform a QR decomposition of \(L^{n+1}_i\) to obtain \(V^{n+1}_j\) and \(S^{n+1}_{ij}\).

\textbf{Cost:} $\bO(r^{2}N_{v})$.
\end{enumerate}

Adding together all of our computational complexity estimates, we get a total cost of $\bO(r^{4} + r^{2}N_{x} + r^{2}N_{v})$ --- compare this with the cost of the full tensor method $\bO(N_xN_v)$.

\section{Fully discrete algorithm}

In this section we address the question of physical and velocity space discretization. One of the
virtues of the dynamical low-rank method is that it decouples the 
discretization of the two bases, $X_{i}$ and $V_{j}$, which may be treated more
or less independently. The $X$ basis is updated by solving a system of coupled
hyperbolic PDEs in (\ref{eqn:K_advance}), while the $V$ basis is updated by
solving a parabolic system in (\ref{eqn:L_IMEX}). These systems are coupled via
the matrix of singular values $S$, as well as weighted inner products of
whichever basis is being held constant (viz. (\ref{eqn:defn_c_ints}),
(\ref{eqn:defn_d_ints}), etc.) We are free to choose whichever discretization is
most appropriate for the corresponding evolution equation of each basis. In this
work we use second-order finite difference discretizations in both $x$ and $v$
for simplicity. In principle, it is easy to choose, for example, a Fourier
spectral method to take advantage of periodicity in the $x$ direction, or even a
more involved method such as Discontinuous Galerkin along one or the other
basis, without increasing the implementation complexity too greatly.

\label{sec:org30a962b}
\subsection{Spatial discretization}
\label{sec:org9fa2cd7}
Our spatial discretization in $x$ is designed to solve the explicit part of the
evolution equation for $K$, which is (\ref{eqn:K_advance}). This is a linear
hyperbolic PDE with the flux matrix $\mathbf{c}^{1}_{jl}$, which is a symmetric matrix.
We opt for a second-order finite difference discretization with flux limiting, as described in \cite{levequeNumericalMethodsConservation1992}, section 16.2.
To illustrate, we consider the situation in two spatial dimensions, $d=2$.
The matrix \(\mathbf{c}^1_{jl}\) consists of
components \(\mathbf{c}^{1;m}_{jl}\) for \(m \in \{1, 2\}\), acting on the \(x\)
and \(y\) directions respectively. The matrices \(\mathbf{c}^{1;m}_{jl}\) are symmetric and real;
recall their definition (\ref{eqn:defn_c_ints}). Therefore they are unitarily
diagonalizable, and we can write
\begin{equation}
(T^m)^{T} \mathbf{c}^{1;m} T^m = \sum_{jl} T^m_{ij} \mathbf{c}^{1;m}_{jl} T^m_{kl} = \lambda_i^m \delta_{ik} = \Lambda^m.
\end{equation}
Left-multiplying (\ref{eqn:K_advance})
by \((T^1)^{T}\), and introducing the eigenbasis \(\hat{K}^n_i = [(T^1)^{T} K^n]_i\), gives the system
\begin{align}
\sum_{jl} T^1_{ij} \left[ I - \frac{\Delta t}{\epsilon} \mathbf{A}^2_{jl} \right] K^{n+1}_l = \hat{K}^n_i - \Delta t \lambda^1_i \partial_x \hat{K}^n_i - \Delta t \sum_{jl} T^1_{ij} \left( \mathbf{c}^{1;2}_{jl} \partial_y K^n_l + \mathbf{A}^1_{jl} K^n_l \right).
\end{align}
At a grid point $x_{p}$, the flux-limited finite difference discretization approximates $\lambda_{i}^{1} \partial_{x} \hat{K}^{n}_{i}(x_{p})$
by a difference of fluxes at half grid points $x_{p+1/2}, x_{p-1/2}$:
\begin{equation}
  \lambda_{i}^{1} \partial_{x} \hat{K}^{n}_{i}(x_{p}) \approx \frac{F(\lambda_{i}, \hat{K}^{n}_{i})_{p+1/2} - F(\lambda_{i}, \hat{K}^{n}_{i})_{p-1/2}}{\Delta x}.
\end{equation}
The flux $F$ is given by the combination of a first-order flux (upwinding) flux,
\begin{equation}
  F_{L}(\lambda, \hat{K})_{p+1/2} = \frac{\lambda}{2}(\hat{K}_{p+1} + \hat{K}_{p}) - \frac{|\lambda|}{2}(\hat{K}_{p+1} - \hat{K}_{p}),
\end{equation}
with a second-order Lax-Wendroff flux. The combination is governed by a
flux-limiter $\phi(\theta)$, which stabilizes the scheme in the presence of sharp changes in the gradient:
\begin{equation}
  F(\lambda, \hat{K})_{p+1/2} = F_{L}(\lambda, \hat{K})_{p+1/2} + \frac{1}{2} \phi(\theta_{p+1/2}) \left( \text{sgn}(\lambda) - \frac{\lambda \Delta t}{\Delta x} \right) \lambda \delta(\hat{K})_{p+1/2},
\end{equation}
where $\delta(\hat{K})_{p+1/2} = (\hat{K}_{p+1} - \hat{K}_{p})$
The quantity $\theta_{p+1/2}$ measures how quickly the gradient is changing in the vicinity of $x_{p+1/2}$, and itself uses upwinding based on the sign of $\lambda$:
\begin{equation}
  \theta_{p+1/2} = \frac{\delta(\hat{K})_{p+1/2 - \text{sgn}(\lambda)}}{\delta(\hat{K})_{p+1/2}}.
\end{equation}
The function $\phi: \mathbb{R} \rightarrow [0, 2]$ is called the limiter, and
there are many options to choose from. We use the Van Leer limiter,
\begin{equation}
  \phi(\theta) = \frac{|\theta| + \theta}{1 + |\theta|}.
\end{equation}

After approximating the term $\lambda_{i}^{1} \partial_{x} \hat{K}^{n}_{i}$, for
each eigenvalue $\lambda_{i}$, we can transform back to the original variables
by left-multiplying with $(T^{1})$:
\begin{align}
\label{eqn:K_spatial_original_variables}
\sum_l \left[ I - \frac{\Delta t}{\epsilon} \mathbf{A}^2_{jl} \right] K^{n+1}_l = K^n_j - \Delta t \sum_i T^1_{ij} \delta_x(\lambda_i, \hat{K}^n)_i - \Delta t \sum_l \left( \mathbf{c}^{1;2}_{jl} \partial_y K^n_l + \mathbf{A}^1_{jl} K^n_l \right).
\end{align}
The discretization in \(y\) is handled similarly, by left-multiplying
(\ref{eqn:K_spatial_original_variables}) by \(T^2\). The above scheme is
second-order in smooth regions of the solution, and degrades to first order
around discontinuities and extrema.

\subsection{Velocity discretization}
\label{sec:org17f3e20}
Our discretization in $v$ is designed to effectively solve (\ref{eqn:L_IMEX}), which is a parabolic system (strictly speaking, a convection-diffusion type equation).
We recall the linear system to be solved here:
\begin{equation}
\label{eqn:L_system}
\sum_{k} \left[ \delta_{ik} - \frac{\Delta t}{\epsilon} (\delta_{ik} \nabla_v^2 - \delta_{ik} v \cdot \nabla_{v} + \mathbf{e}^\star_{ik} \cdot \nabla_v) \right]L_{k} = RHS.
\end{equation}
It is convenient to discretize this operator by splitting the left hand side into a diagonal ($i=k$) term and an offdiagonal term.
In the case when $i = k$, we have
\begin{align}
  \left(\delta_{ik} \nabla_{v}^{2} - \delta_{ik} v \cdot \nabla_{v} + \mathbf{e}^{\star}_{kk} \cdot \nabla_{v} \right) L_{k} &= \left(\nabla_{v}^{2} - v\cdot \nabla_{v} + \mathbf{e}^{\star}_{kk}\cdot\nabla_{v}\right)L_{k} \nonumber \\
  &= \left[((\nabla_{v} - v) + \mathbf{e}^{\star}_{kk})\cdot\nabla_{v} \right] L_{k} \nonumber \\
  &= \frac{1}{M^{k}}\nabla_{v} \cdot (M^{k}\nabla_{v} L_{k}) \nonumber \\
  \label{eqn:L_diagonal}
  &:= \mathbb{T}_{k}(\mathbf{e}^{\star}_{kk}) L_{k},
\end{align}
where $M^{k}$ is the local Maxwellian
\begin{equation}
  M^{k} = e^{-\frac{|v - \mathbf{e}^{\star}_{kk}|^{2}}{2}}.
\end{equation}
In one dimension, a second-order-accurate central difference discretization of
(\ref{eqn:L_diagonal}) is
\begin{align}
  (\mathbb{T}_{k}(\mathbf{e}^{\star}_{kk}) L)_{p} &\approx \frac{1}{M^{k}_{p}\Delta v} \left( M^{k}_{p+1/2} \frac{L_{p+1}-L_{p}}{\Delta v} - M^{k}_{p-1/2} \frac{L_{p-1}-L_{p}}{\Delta v} \right) \\
  &= \frac{M^{k}_{p+1/2} L_{p+1} - (M^{k}_{p+1/2} + M^{k}_{p-1/2})L_{p} + M^{k}_{p-1/2} L_{p-1}}{M^{k}_{p}\Delta v^{2}}.
\end{align}
The off-diagonal terms are simply
\begin{align}
  \mathbf{e}_{ik}^{*} \cdot \nabla_{v} L_{k} := \mathbb{U}_{ik} L_{k}.
\end{align}
For simplicity we discretize this using a second-order centered difference
operator. Stability is not a concern, since it will be coupled to a Backwards
Euler timestepping scheme. With these discretizations in hand the implicit step
for $L$ takes the form
\begin{equation}
  \sum_{k} \left[ \delta_{ik} - \frac{\Delta t}{\epsilon} (\delta_{ik} \mathbb{T}_{k}(\mathbf{e}^{\star}_{kk}) + \mathbb{U}_{ik})\right] L_{k}^{n+1} = L_{i}^{n} - \Delta t \sum_{k}(v \cdot \mathbf{d}_{ik}^{\star} + \hat{c}_{ik}) L_{k}^{n}.
\end{equation}
Despite being of size $rN_{v} \times rN_{v}$, this linear system is quite
sparse, having roughly $\bO(r^{2}N_{v})$ nonzero entries. It is therefore amenable
to fast solution by iterative solvers. Since it is not symmetric, we use the
Restarted GMRES \cite{saadGMRESGeneralizedMinimal1986} iterative algorithm. We
find good results by preconditioning with the constant matrix
$\left( \delta_{ik} - \delta_{ik}\frac{\Delta t}{\epsilon}\mathbb{T}_{k}(\mathbf{0})\right)^{-1}$.
Timings of our code indicate that this step takes on the same order of magnitude
as the other components of the algorithm, up to the largest problems we consider
here.

\section{Asymptotic behavior of the discrete scheme}

In this section we demonstrate that the discrete scheme described in the preceding
sections preserves the asymptotic limit (\ref{eqn:limiting_system}) as $\epsilon \rightarrow 0$.

We consider the limit of the discrete system at the level of $g$, which is advanced via the $K$, $S$, and $L$ steps
with an accuracy that is first-order in time and second-order in space: 
\begin{equation} \label{scheme}
  \frac{g^{n+1} - g^{n}}{\Delta t} = -v \cdot \nabla_{x} g^{n} - \mathcal{M}^{n}g^{n} + \frac{1}{\epsilon} (M^{n})^{-1}\nabla_{v}\cdot(M^{n}\nabla_{v} g^{n+1}) + \bO(\Delta t + \Delta x^{2}).
\end{equation}
Furthermore, the electric field is advanced by
\begin{equation} \label{EE}
\frac{E^{n+1}-E^n}{\Delta t}=-J^n.
\end{equation}

From (\ref{scheme}), we can see that
\begin{equation} \label{gg}
\bO(\epsilon)=(M^{n})^{-1}\nabla_{v}\cdot(M^{n}\nabla_{v} g^{n+1}) \ \Longrightarrow \  g^{n+1}=c+\bO(\epsilon),
\end{equation}
that is to say, after one time step, we expect the solution $g$ to be close to a constant function in $v$ when $\epsilon$ is small. To see this, just note the following
\begin{equation}
\int \nabla_v \cdot  [M \nabla_v g] \log g \rd{v}=\int \nabla_v  \cdot [Mg\nabla_v  \log g]\log g \rd{v}=-\int  Mg \left| \nabla_v \log g\right|^2\rd{v}\leq 0,
\end{equation}
where the equality holds if and only if $g$ is a function independent of $v$. Moreover, (\ref{gg}) implies 
\begin{equation}
\rho^{n+1}=\int M^{n+1}g^{n+1}\rd{v}=c\int M^{n+1}\rd{v}+\bO(\epsilon)=c+\bO(\epsilon),
\end{equation}
\begin{equation}
  J^{n+1} =\int vM^{n+1}g^{n+1}\rd{v}=c \int v M^{n+1}\rd{v}+\bO(\epsilon)=cE^{n+1}+\bO(\epsilon)=\rho^{n+1}E^{n+1} + \bO(\epsilon).
\end{equation}

On the other hand, we can multiply (\ref{scheme}) by $M^{n}$ and integrate in $v$ to obtain
\begin{align}
  \frac{\int M^{n} g^{n+1}\rd{v} - \rho^{n}}{\Delta t} &= -\int v \cdot (\nabla_{x} g^{n}) M^{n} \rd{v} - \int (\partial_{t} M^{n} + v \cdot \nabla_{x} M^{n})g^{n} \rd{v} \nonumber \\
                                                     &\qquad+ \frac{1}{\epsilon} \cancel{\int \nabla_{v} \cdot (M^{n} \nabla_{v} g^{n+1})\rd{v}} + \bO(\Delta t + \Delta x^{2}) \nonumber \\
                                                     &= -\int [v \cdot (\nabla_{x} g^{n}) M^{n} + ( v \cdot \nabla_{x}M^{n} )g^{n}] \rd{v} - \partial_{t} E^n \cdot \int (v-E^n) M^{n} g^{n}\rd{v} \nonumber \\
                                                     &\qquad+ \bO(\Delta t + \Delta x^{2}) \nonumber \\
  &= -\nabla_{x}\cdot J^{n} - c\partial_{t} E^n \cdot \cancel{\int (v-E^n) M^{n} \rd{v}}+\bO(\epsilon+\Delta t + \Delta x^{2}),  \label{rhorho}
\end{align} 
where we used $g^n=c+\bO(\epsilon)$ for $n\geq 1$.

Finally, noticing that
\begin{equation}
  \rho^{n+1} - \int M^{n}g^{n+1}\rd{v} = \int (M^{n+1} - M^{n})g^{n+1}\rd{v} = c\int (M^{n+1} - M^{n})\rd{v} +\bO(\epsilon)=\bO(\epsilon),
\end{equation}
and using $J^n=\rho^n E^n+\bO(\epsilon)$ for $n\geq 1$ (\ref{rhorho}) becomes
\begin{equation}
  \label{eqn:rho_asymptote}
  \frac{\rho^{n+1} - \rho^{n}}{\Delta t} = -\nabla_{x}\cdot (\rho^nE^n) + \bO\left(\epsilon + \Delta t + \Delta x^{2}+\frac{\epsilon}{\Delta t}\right).
\end{equation}
(\ref{EE}) becomes
\begin{equation}
\label{eqn:E_asymptote}
\frac{E^{n+1}-E^n}{\Delta t}=-\rho^nE^n+\bO(\epsilon).
\end{equation}
Equations (\ref{eqn:rho_asymptote}) and (\ref{eqn:E_asymptote}) form a first-order in time discretization of
(\ref{eqn:limiting_system}) as $\epsilon \rightarrow 0$, as expected.

\section{Numerical results}
In this section, we present extensive numerical results in 1D1V and 2D2V to illustrate the accuracy and efficiency of the proposed low-rank algorithm. We will see that the algorithm becomes quite low rank in the asymptotic limit $\epsilon \rightarrow 0$ by our design. On the other hand, in the kinetic and transition regimes, the numerical rank needed appears higher but still relatively small compared to $N_x$ or $N_v$. Therefore, the proposed algorithm presents as a very effective method for the Vlasov-Amp\`{e}re-Fokker-Planck system over a wide range of problems.

\input{1d1v}
\input{2d2v}

\section{Conclusion}
We have proposed and implemented an efficient algorithm for the electrostatic
Vlasov equation with linear Fokker-Planck collision operator. By dividing by the
Maxwellian, we are able to represent the quotient with a low-rank approximation,
thereby capturing the fluid limit with very little computational effort.
Moreover, our method is also efficient when far from the fluid limit, owing to
the great reduction in computational complexity afforded by the dynamical
low-rank method. In order to get an efficient overall algorithm, we used the
fact that the Maxwellian limit of our equation is isothermal to quickly compute
a convolution with the Fast Fourier Transform. Our implementation is found to be
multiple orders of magnitude faster than a full-tensor numerical solution, with
better asymptotic scaling and constant coefficients for moderately sized problems.

\bibliographystyle{plain}
\bibliography{references-bibtex.bib}

\appendix
\section{Appendix: Timestepping for the backwards-in-time $S$ step}
\label{S_STEP}
\input{s_step}

\end{document}

%% file: 1d1v.tex
\subsection{1D1V examples}

\subsubsection{Convergence study}

\begin{figure}
\begin{subfigure}{.48\textwidth}
  \centering
  \includegraphics[width=\textwidth]{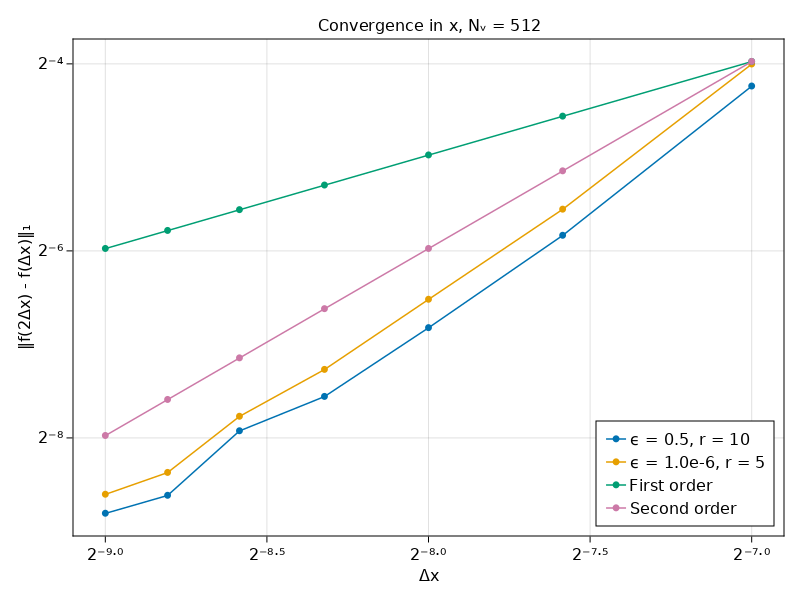}
  \caption{Second order convergence in $x$. $N_{v} = 512$, $N_{x}$ between 64 and 512.}
  \label{fig:x_convergence}
\end{subfigure}
\hfill
\begin{subfigure}{.48\textwidth}
  \centering
  \includegraphics[width=\textwidth]{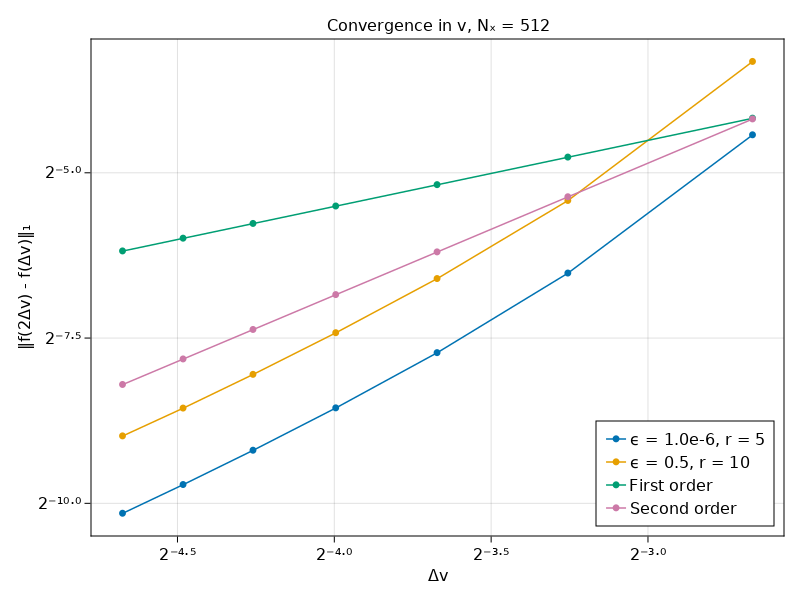}
  \caption{Second order convergence in $v$. $N_{x} = 512$, $N_{v}$ between 64 and 512.}
  \label{fig:v_convergence}
\end{subfigure}
\caption{Convergence of the solution as a function of the grid spacing in both $x$ and $v$. We achieve second-order convergence in both coordinates.}
\end{figure}

To verify second order convergence of our scheme in physical space and velocity
space, we perform convergence studies comparing the relative errors in $f$ as
the grid is refined. In order to evaluate the performance of the discrete Fokker-Planck
collision operator, we use a nonequilibrium initial condition. The initial
distribution consists of two counterstreaming beams moving at velocities
$\pm 1.5$.
\begin{align}
  \label{eqn:counterstreaming}
  &f(x, v, 0) = \frac{\rho_{0}(x)}{2\sqrt{2\pi}} \left[ e^{\frac{-|v - 1.5|^{2}}{2}} + e^{\frac{-|v + 1.5|^{2}}{2}} \right], \\
  &\rho_{0}(x) = \sqrt{2\pi} (2 + \cos(2\pi x)).
\end{align}
The initial electric field is determined by Poisson's equation,
\begin{align}
  \label{eqn:poissons}
  &E_{0} = -\nabla_{x} \phi(x), \quad -\nabla_{x}^{2} \phi(x) = \rho_{0}(x) - \eta(x), \\
  &\eta(x) = \frac{2\sqrt{2\pi}}{1.2661} e^{\cos(2\pi x)}.
\end{align}
We evaluate the convergence in both the kinetic and fluid
regimes, with $\epsilon = 0.5$ and $\epsilon = 10^{-6}$, respectively. The fluid
regime is adequately resolved with $r = 5$, while the kinetic regime requires a
higher rank of $r = 10$. The spatial domain is periodic on the interval $[0, 1]$,
and the velocity domain is the interval $[-10, 10]$.
Convergence is verified by holding one of $N_{x}, N_{v}$ fixed, while the
other is varied.
The timestep $\Delta t$ is chosen to give a CFL number of $0.25$ at the finest grid,
for which we use $\Delta t = \frac{1.0}{(512)(4v_{max})} = \num{4.88e-5}$. This is
found to be sufficient for the spatial discretization error to dominate.
The $L^{1}$ norm of the difference between $f$ at subsequent
levels of approximation is computed by linearly interpolating the solution at
the coarser grid onto the finer grid. The successive differences are observed to
converge at second order in the grid spacing, namely $\Delta x$
(Figure \ref{fig:x_convergence}) or $\Delta v$ (Figure \ref{fig:v_convergence}).

\subsubsection{Verification that the asymptotic limit is rank-1}
\begin{figure}
\begin{subfigure}{.5\textwidth}
  \centering
  \includegraphics[width=\linewidth]{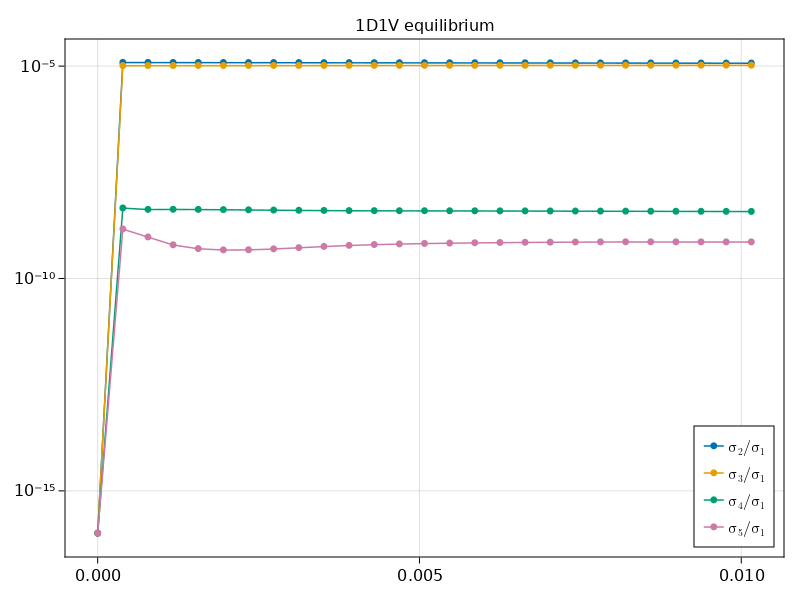}
  \caption{1D1V initial condition in local equilibrium}
  \label{fig:rank_history_local_eq}
\end{subfigure}
\begin{subfigure}{.5\textwidth}
  \centering
  \includegraphics[width=\linewidth]{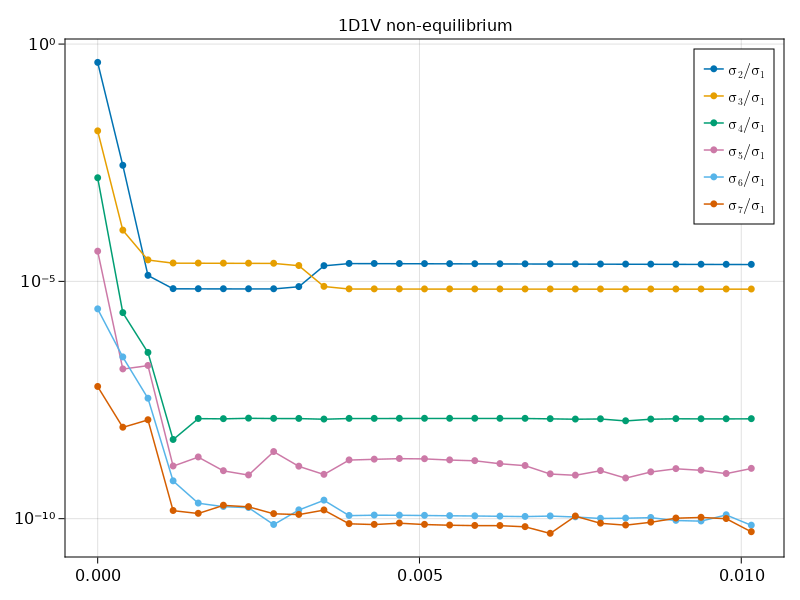}
  \caption{1D1V initial condition in non-equilibrium}
  \label{fig:rank_history_counterstreaming}
\end{subfigure}
\caption{Evolution of normalized singular values of the solution in the fluid regime, $\epsilon = 10^{-6}$. There is a clear separation, of five orders of magnitude, between the largest and second largest singular values. This holds whether the solution begins in local equilibrium (left) or in the non-equilibrium initial condition (\ref{eqn:counterstreaming}).}
\end{figure}

In order to verify that our method efficiently captures the asymptotic limit
(\ref{eqn:limiting_system}), we examine the evolution of the singular values of
the low rank solution for very small $\epsilon$ ($\epsilon = 10^{-6}$). The
singular values of the low rank solution are simply the diagonal entries of $S$.
We consider a solution beginning in local equilibrium,
\begin{align}
  &f(x, v, 0) = \frac{\rho_{0}(x)}{\sqrt{2\pi}} e^{-\frac{|v-E_{0}|^{2}}{2}}, \\
  &\rho_{0}(x) = \frac{\sqrt{2\pi}}{2} (2 + \cos(2\pi x)),
\end{align}
where $E_{0}$ satisfies (\ref{eqn:poissons}) with
\begin{align}
  \eta(x) = \frac{\sqrt{2\pi}}{1.2661} e^{\cos(2\pi x)}.
\end{align}
We evolve the initial condition with a rank of 5, until time $t = 0.01$, which
is enough to demonstrate that the asymptotic limit is captured.
$N_{x} = N_{v} = 128$ grid points are used in each direction. The timestep
chosen is $\Delta t = \num{3.9e-4}$. The evolution of singular values in Figure
\ref{fig:rank_history_local_eq} shows that the solution maintains a clear
separation between the first singular value and the rest. Figure
\ref{fig:rank_history_counterstreaming} demonstrates the same behavior, but for
a solution beginning in the counterstreaming beams initial condition,
(\ref{eqn:counterstreaming}), evolved with rank 10. The solution takes slightly
longer to ``settle down'', but after three time steps it shows the same rank
separation as in the equilibrium case.

\subsubsection{Comparison of fluid and kinetic regimes}

Conversely to the clear singular value separation observed in the fluid regime for
very small $\epsilon$, solutions in the kinetic regime exhibit slower singular
value decay. To demonstrate, we consider a ``bump-on-tail'' initial condition
evolving in both the fluid and kinetic regimes. In the fluid regime, the bump
disappears within a single timestep, and the slow time scale dynamics of the
limiting fluid equation take over. In the kinetic regime, we observe the
shearing behavior characteristic of low-collision phase space flows. Our method
is designed to capture the low-rank structure inherent in the fluid equation,
and so it is not surprising that the kinetic solution requires a higher rank to
capture effectively.

The ``bump-on-tail'' distribution is defined by
\begin{equation}
  f_{0}(x, v) = \frac{\rho_{0}(x)}{(2\pi)^{1/2}} \left( e^{-\frac{|v|^{2}}{2}} + e^{-\frac{|v-1.5|^{2}}{2 T_{\text{cold}}}} \right),
\end{equation}
where the temperature of the perturbation is $T_{\text{cold}} = .005$.
The density $\rho_{0}(x)$ is initialized to a
Gaussian pulse centered at $x = 0.3$. To induce dynamics in the limiting fluid
equation, we initialize the background charge density $\eta$ with a potential
well centered at $x = 0.6$:
\begin{equation}
  \rho_{0}(x) = 0.3 + e^{-\frac{|x-0.3|^{2}}{0.01}}, \quad \eta(x) = 0.3 + e^{-\frac{|x-0.6|^{2}}{0.01}}.
\end{equation}
The electric field is initialized via the solution to Poisson's equation
(\ref{eqn:poissons}), as above.

\begin{figure}
\begin{subfigure}{0.5\textwidth}
  \centering
  \includegraphics[width=\textwidth]{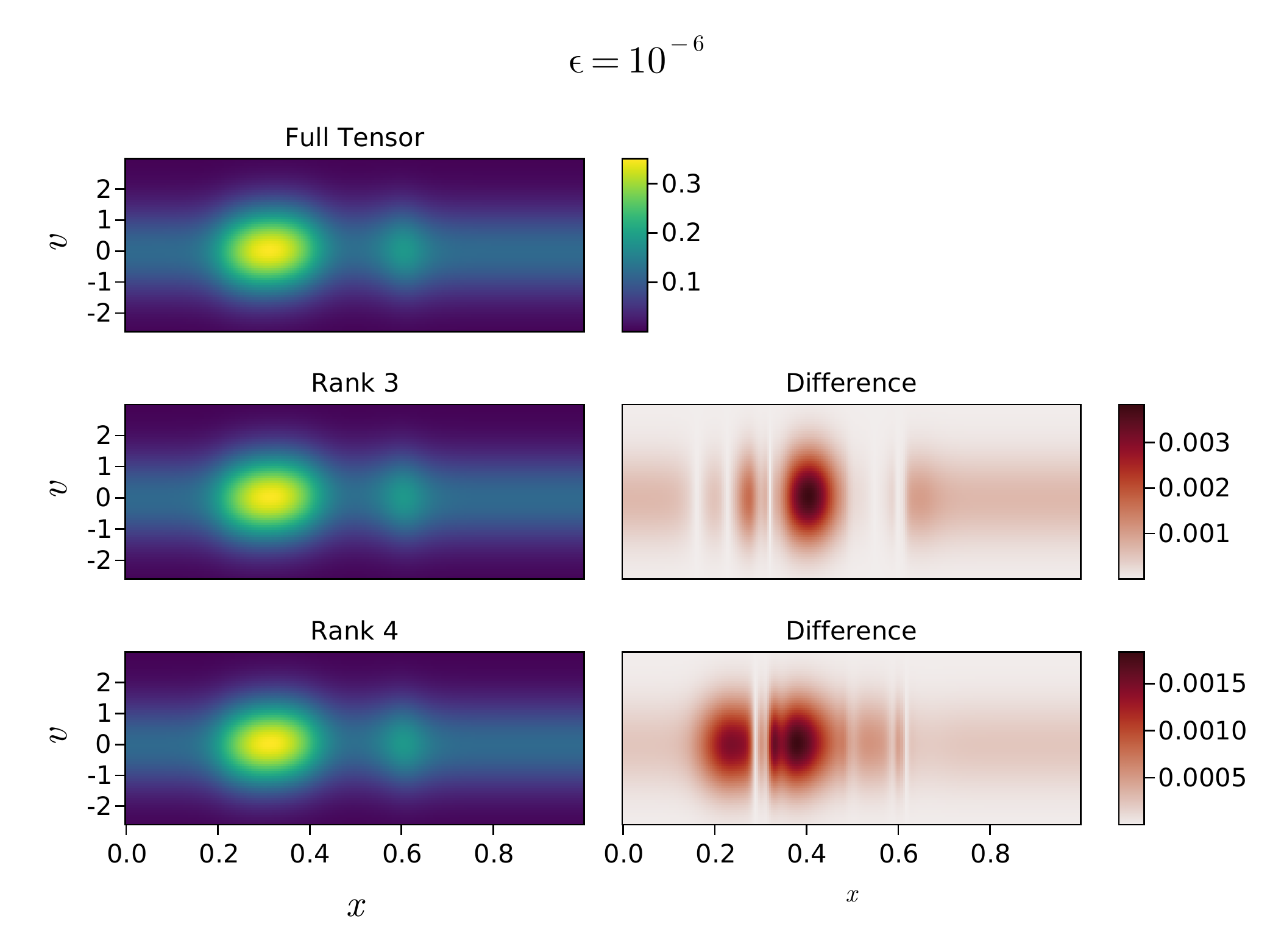}
  \caption{Fluid regime solutions at ranks $r = 3, 4$. Good accuracy is achieved with only a handful of ranks.}
  \label{fig:rank_comparison_fluid}
\end{subfigure}
\begin{subfigure}{0.5\textwidth}
  \centering
  \includegraphics[width=\textwidth]{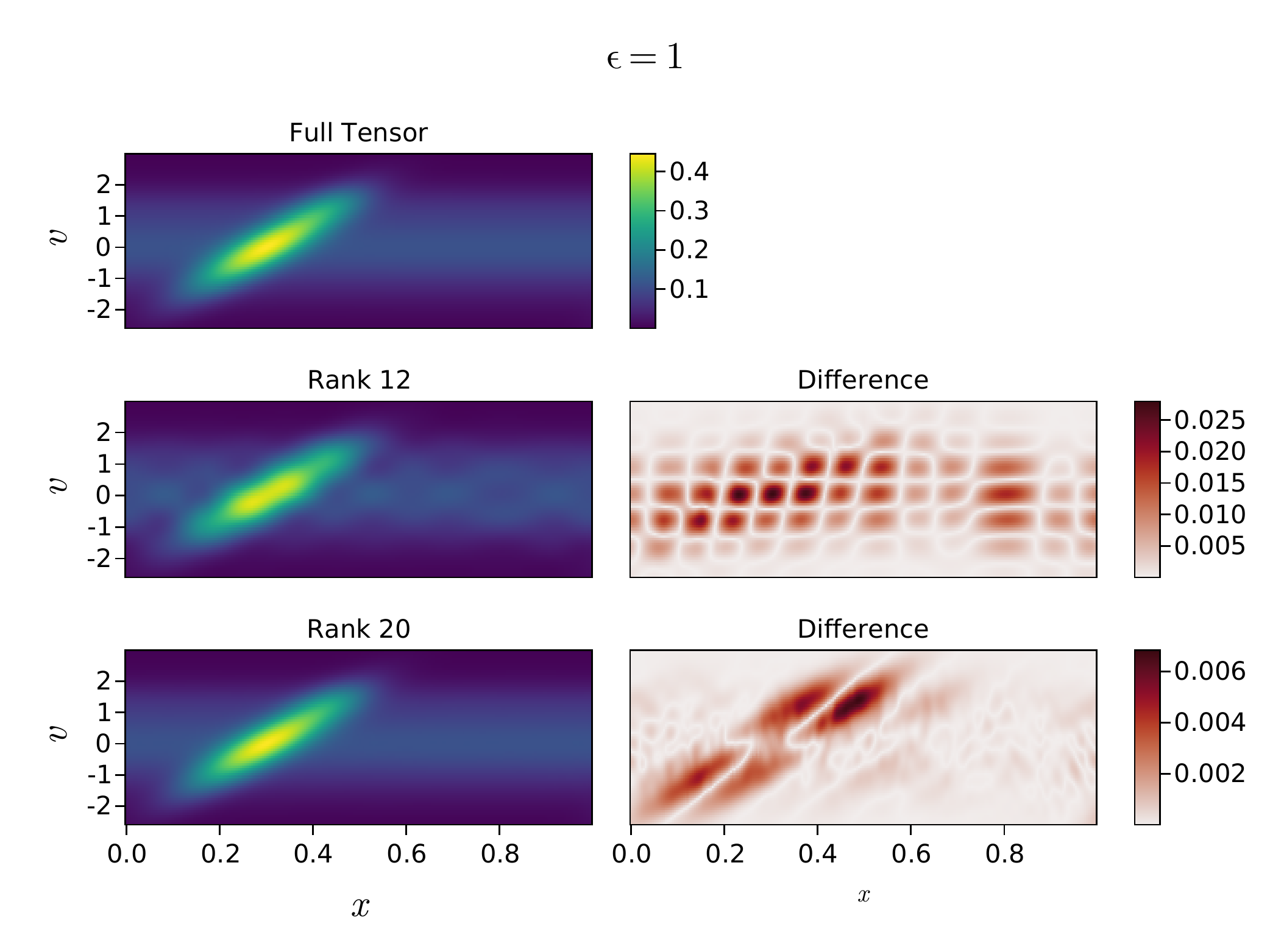}
  \caption{Kinetic regime solutions. Significantly more ranks are required to evolve the solution far from local equilibrium.}
  \label{fig:rank_comparison_kinetic}
\end{subfigure}
\caption{Comparison of the rank $r$ required for a solution in the fluid (left) and kinetic (right) regime. The fluid solution with $\epsilon = 10^{-6}$ requires significantly lower rank, demonstrating the ability of our approach to capture the fluid limit without excessive computation, by seeking a low-rank approximation to $g = f/M$ instead of $f$ itself. Solutions are computed with $N_{v} = 128$ and $N_{x} = 256$. }
\label{fig:rank_comparison}
\end{figure}

The solution is computed with $N_{x}=256$ grid points in the $x$ coordinate, and
$N_{v} = 128$ in the $v$ coordinate. We use $\epsilon=10^{-6}$ to demonstrate
the fluid regime behavior, and in the kinetic regime we use $\epsilon=1$. Both
of our low-rank solutions are compared to a full tensor solution of
(\ref{eqn:vpfp_combined}) with the same discretization parameters, using the
scheme proposed in
\cite{jinAsymptoticPreservingScheme2011}. The solution is run until $t = 0.1$ in
the fluid regime, and $t = 0.5$ in the kinetic regime, with time steps of
$\Delta t = \num{3.9e-4}$ in all cases. Results are shown in Figure
\ref{fig:rank_comparison}. We find that as expected, only a handful of ranks are
required to obtain good accuracy in the fluid regime. On the other hand, the
kinetic solution requires around $r=20$ for this problem.

%% file: 2d2v.tex
\subsection{2D2V examples}

\subsubsection{Climbing an electrostatic potential hill}

To demonstrate that our method can handle nontrivial dynamics in the kinetic
regime, we consider the problem of a density ``pulse'' climbing an electrostatic
potential hill. We initialize the background density $\eta$ uniform everywhere
except for a band through the center of the domain where it is set to zero. This
creates a region of negative charge density through which the electron fluid
cannot pass, unless it has enough inertia to do so. Since in the fluid limit
inertial forces vanish, this will only occur in the kinetic regime. We use an
elongated Gaussian initial density centered to the left of the potential hill
and oriented obliquely to the grid:
\begin{equation}
  \rho_{0}(x, t) = 0.1 + \frac{0.0003}{2\pi |\bm{\Sigma}|} e^{-\frac{(x-x_{0})^{T} \bm{\Sigma}^{-1} (x-x_{0})}{2}}, \quad x \in [0, 1]^{2},
\end{equation}
where $x_{0} = [0.3, 0.3]$, $\bm{\Sigma} = \bm{R} \Lambda \bm{R}^{-1}$, $\Lambda$ is a diagonal matrix
with entries $[0.006, 0.03]$, and $\bm{R}$ a rotation matrix through an angle of
$\pi/4$. We initialize a uniform Maxwellian velocity distribution throughout the
domain so that the pulse is traveling along the direction of its major axis:
\begin{equation}
  f_{0}(x, v) = \frac{\rho_{0}(x)}{2\pi T} e^{-\frac{|v - u_{0}|^{2}}{2T}}, \quad v \in [-5, 5]^{2},
\end{equation}
where $T = 0.01$ and $u_{0} = [0.5, 0.5]^{T}$. The background density $\eta$ is initialized constant on its support, which is the entire domain excluding a band between $x_{l} = 0.55$ and $x_{r}=0.7$. This creates the potential hill which a kinetic distribution is able to pass over, while the fluid solution remains on the left side, where it starts.

\begin{figure}
  \begin{subfigure}{0.5\textwidth}
    \centering
    \includegraphics[width=\linewidth]{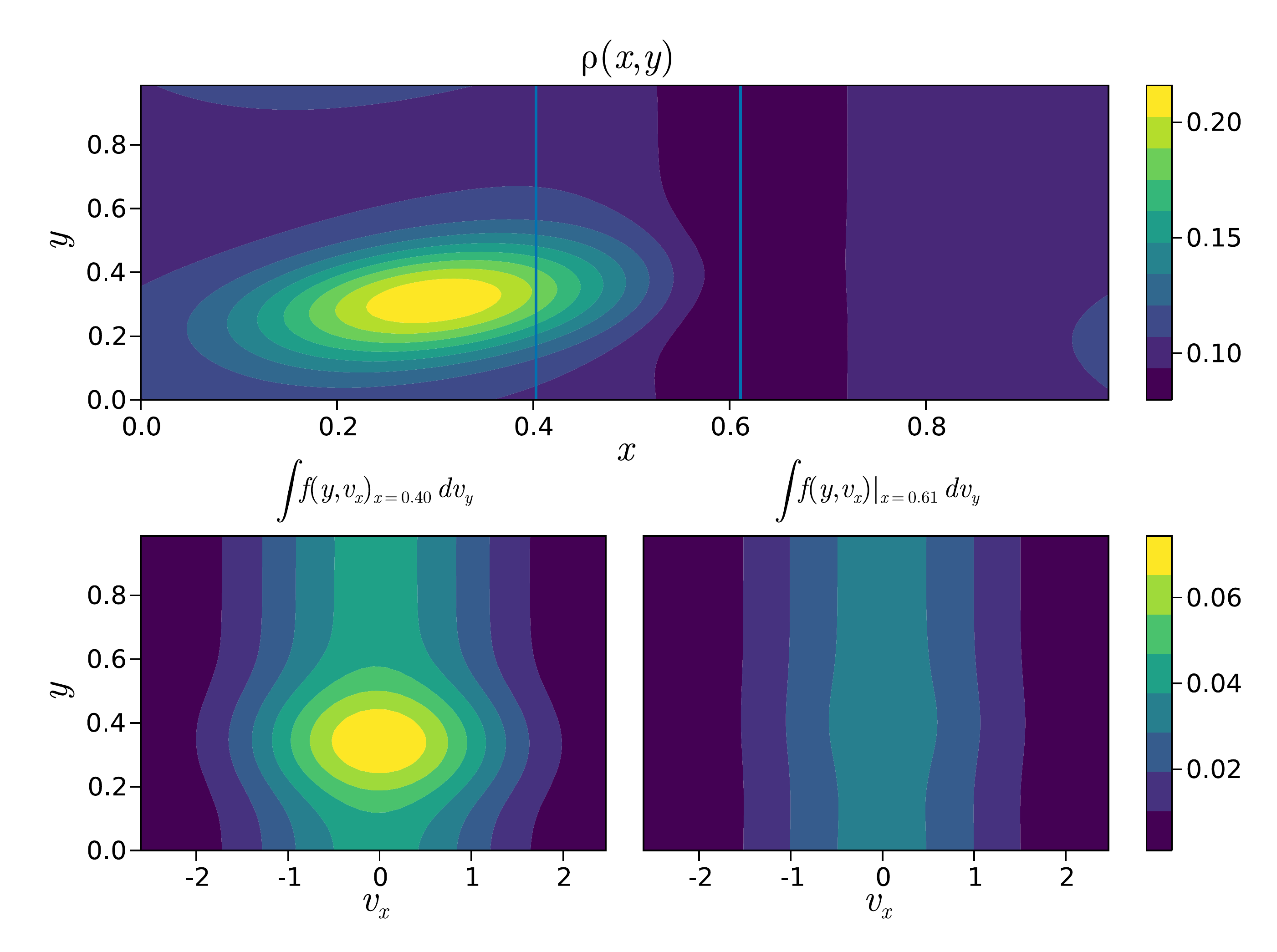}
    \caption{Low rank solver, $\epsilon=0.01, r = 10$}
  \end{subfigure}
  \begin{subfigure}{0.5\textwidth}
    \includegraphics[width=\linewidth]{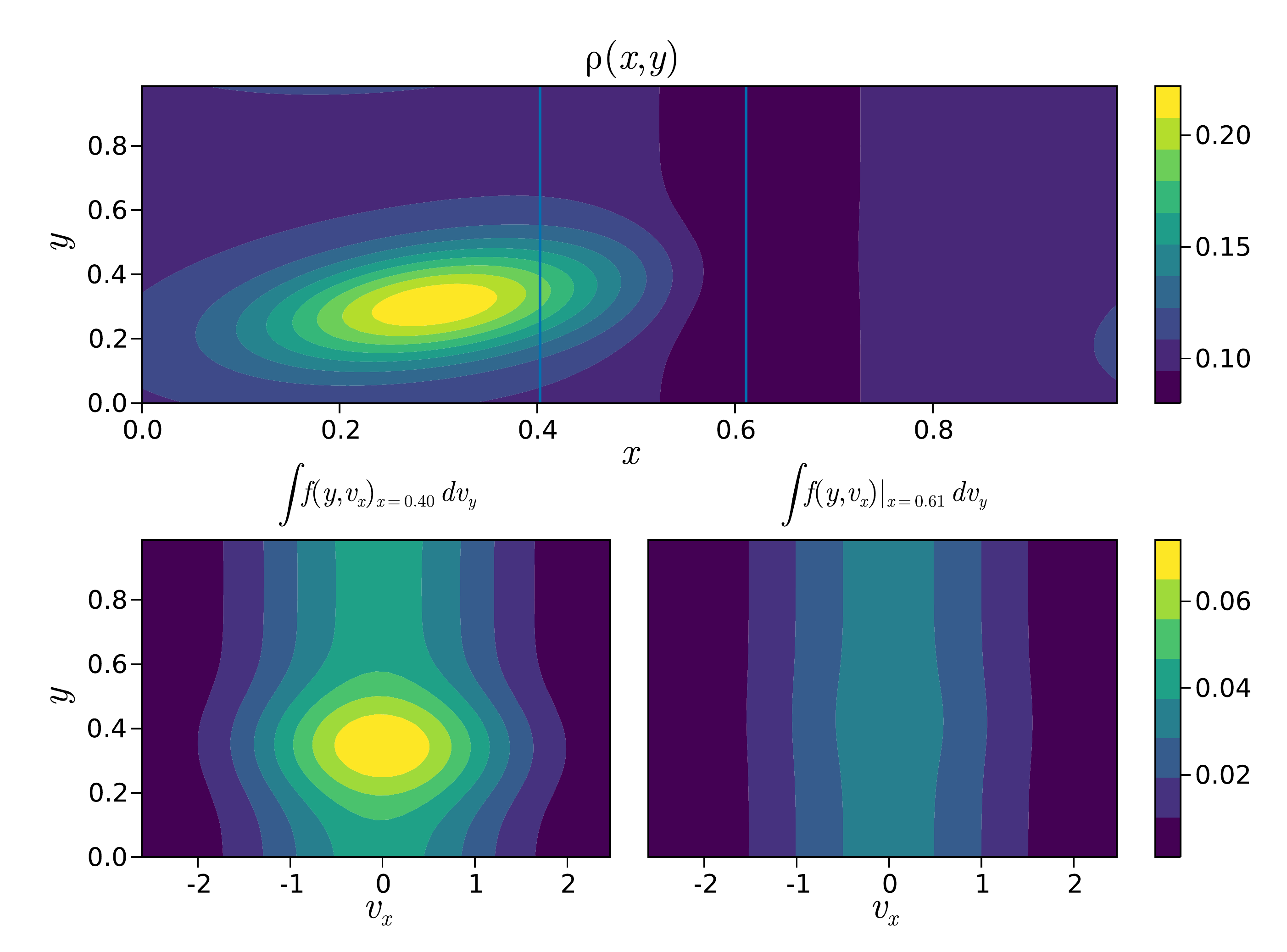}
    \caption{Full tensor solver, $\epsilon=0.01$}
  \end{subfigure}

  \begin{subfigure}{0.5\textwidth}
    \centering
    \includegraphics[width=\linewidth]{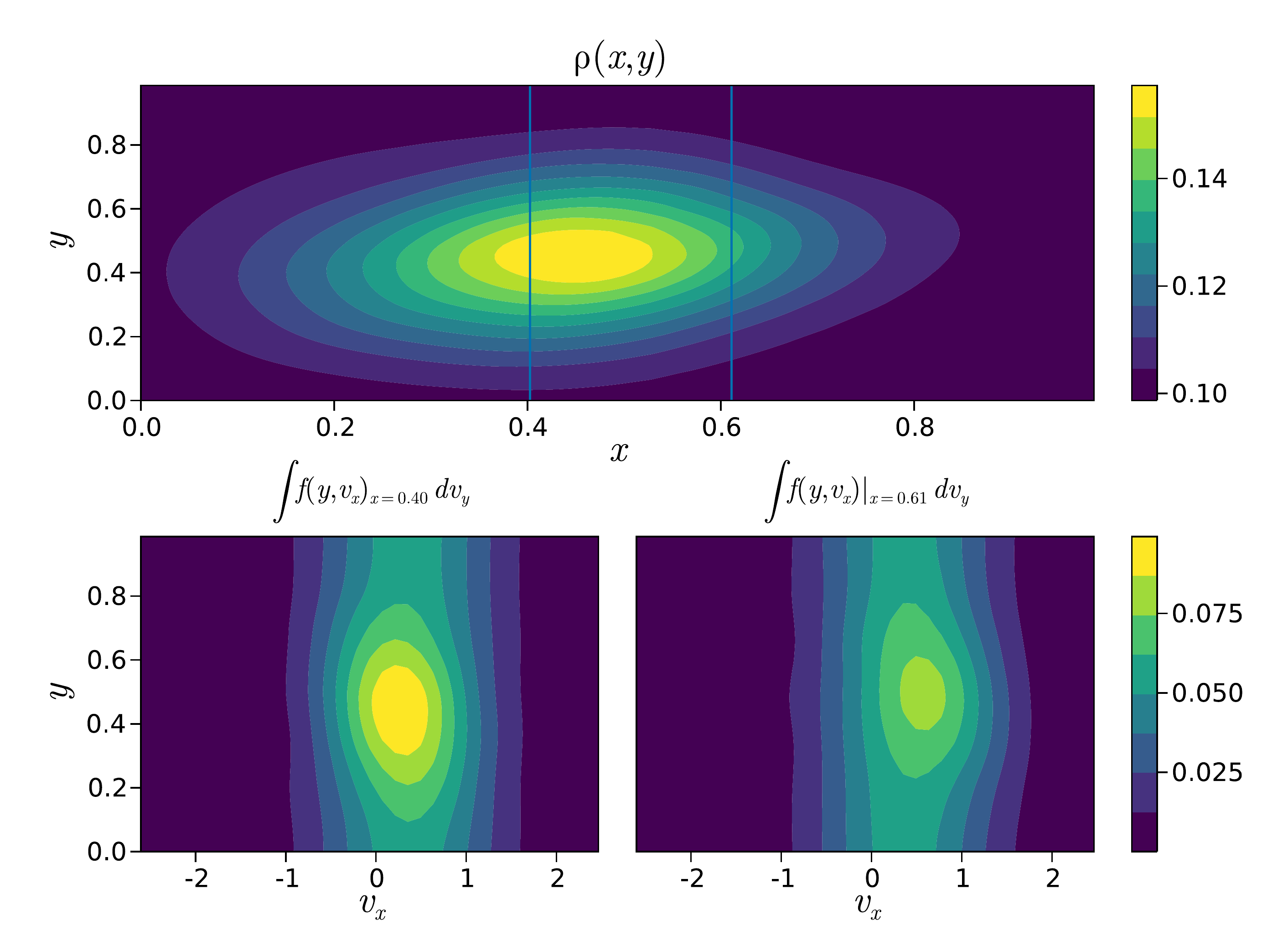}
    \caption{Low rank solver, $\epsilon=1.0, r = 15$}
  \end{subfigure}
  \begin{subfigure}{0.5\textwidth}
    \includegraphics[width=\linewidth]{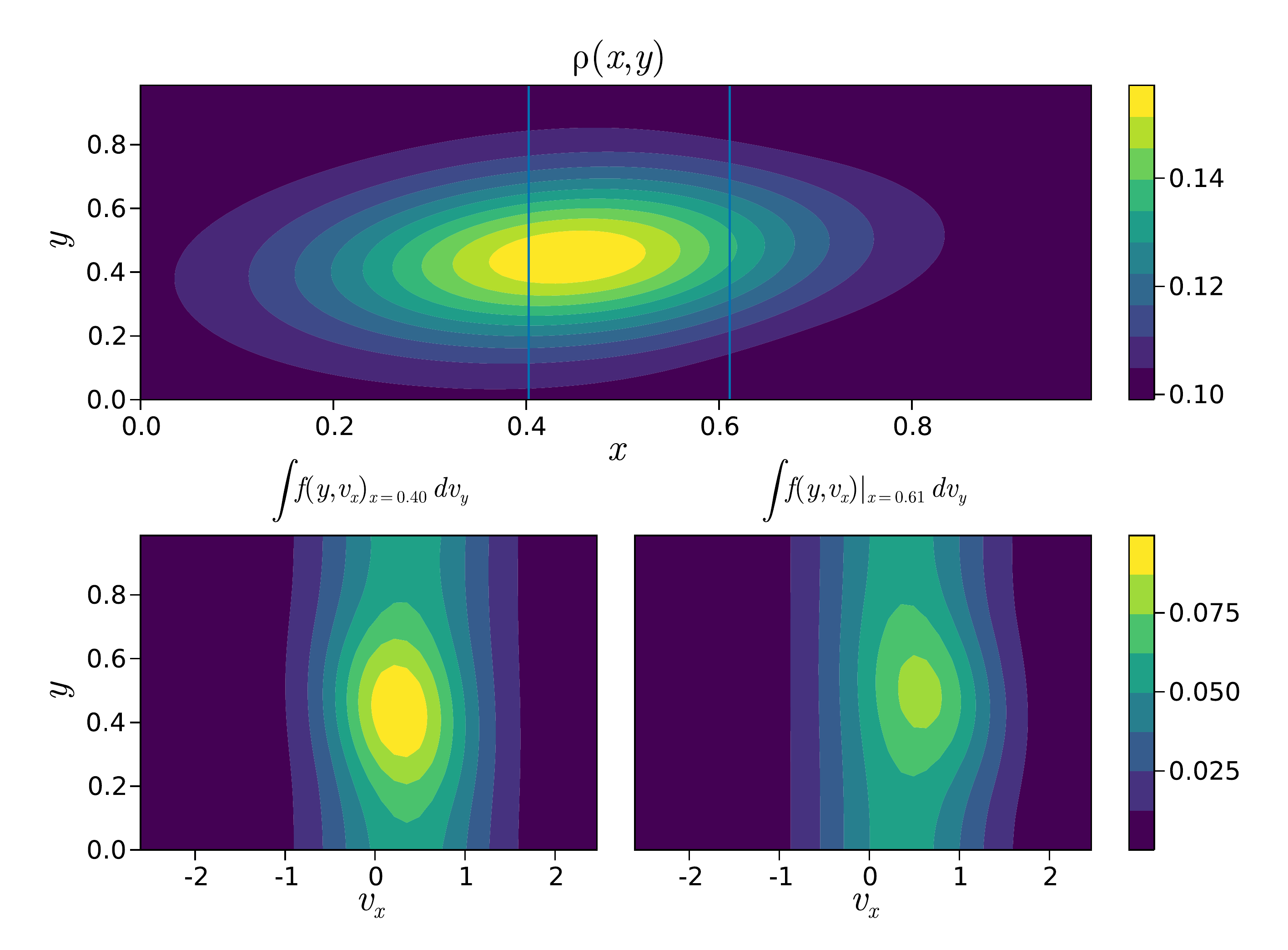}
    \caption{Full tensor solver, $\epsilon=1.0$}
  \end{subfigure}

  \caption{Results of the potential hill problem at $T = 0.35$, comparing low rank (left) and full tensor (right) solutions. Marginal distributions of $f$ as a function of $y, v_{x}$ at $x \approx 0.4$, $x \approx 0.6$ are shown in the bottom pair of plots for each case. We see good agreement with the full tensor solution and with our heuristic predictions for what each flow should do. \label{fig:potential_hill}}
\end{figure}

To illustrate both regimes, we use the values $\epsilon=1.0$ and
$\epsilon=0.01$, and a computational domain with $N=72$ grid points in each of
the four coordinates. We use a fixed time step of $\Delta t = \num{6.9e-4}$. The
results at time $T = 0.35$ are shown in Figure \ref{fig:potential_hill}. As
expected, the kinetic solution retains a significant flow velocity throughout
the domain, and its inertia carries it over the $x = 0.7$ line. On the other
hand the fluid solution is pushed out of the region of negative charge density
by electrostatic forces. The presence of the potential hill in the interval
$x \in [0.55, 0.7]$ is clearly visible in the density plots for the fluid
regime.

Our dynamical low-rank method shows its computational advantages on this 2D2V
problem, even for small problem sizes. Timings for a single timestep are detailed in Table
\ref{table:timings}. We observe that the computational cost of the algorithm
is $\bO(r^{2}N_{x} + r^{2}N_{v})$, compared to the $\bO(N_{x}N_{v})$ of the
full tensor algorithm. The constant factors are small enough to already be dominated
at $N_{x} = N_{v} = 24^{2}$.

\begin{table}
  \centering
  \begin{tabular}{r l l l | l}
    $N$ & Low Rank & & & Full Tensor \\
      & $r=5$ & $r=10$ & $r=15$ \\
    \hline
    24 & 1.00 / -            & 2.14 / -            & 3.89 / -              & 13.8 / - \\
    48 & 3.21 / \textbf{1.7} & 7.94 / \textbf{1.9} & 14.4 / \textbf{1.9}   & 256  / \textbf{4.2} \\
    72 & 10.5 / \textbf{2.1} & 23.6 / \textbf{2.2} & 46.7 / \textbf{2.3}   & \num{1.5e3} / \textbf{4.3} \\
    96 & 18.6 / \textbf{2.1} & 38.3 / \textbf{2.1} & 54.7 / \textbf{1.9}   & \num{5.3e3} / \textbf{4.3} \\
    120 & 29.2 / \textbf{2.1} & 52.2 / \textbf{2.0}  & 88.9 / \textbf{1.9} & \num{1.23e4} / \textbf{4.2}
  \end{tabular}
  \caption{Computational runtime per time step of the 2D2V potential hill problem. $N$ is the number of grid points in each dimension, so that $N_{x} = N_{v} = N^{2}$. Runtimes are normalized to the $r=5, N=24$ size. For example, the runtime of the full tensor step is 13.8 times longer than that of the rank 5 low rank solve when $N = 24$. The bolded numbers are the empirical exponent of $N$. We see agreement with the expected asymptotic complexity of 2 for the low-rank case, compared with 4 for the full tensor solver. \label{table:timings}}
\end{table}

\subsubsection{Relaxation of a cold beam}

To demonstrate the relaxation of the solution towards the local Maxwellian in 2
velocity dimensions, we consider a ``cold beam'' initial condition:
\begin{align}
  &f(x, v, 0) = e^{-\frac{|v - u|^{2}}{0.5}}, \\
  &u = [4, 2]^{T}.
\end{align}
This is discretized on a doubly spatially periodic unit domain
$[0, 1)^2 \times [-10, 10]^2$, with $N_{x}= 32^{2}$ spatial grid points
dimension, and $N_{v} = 128^{2}$ velocity grid points. We choose
$\Delta t \approx 0.0008$, and evolve the distribution until $T = 0.3$. The
difference between the solution $f$ and the local equilibrium distribution $M$
is plotted for three intermediate points in time, along with the history of the
$L^{1}$ norm of the difference. It can be seen that the deviation from local
equilibrium decays exponentially.

\begin{figure}
  \begin{subfigure}{0.8\textwidth}
    \centering
    \includegraphics[width=\linewidth]{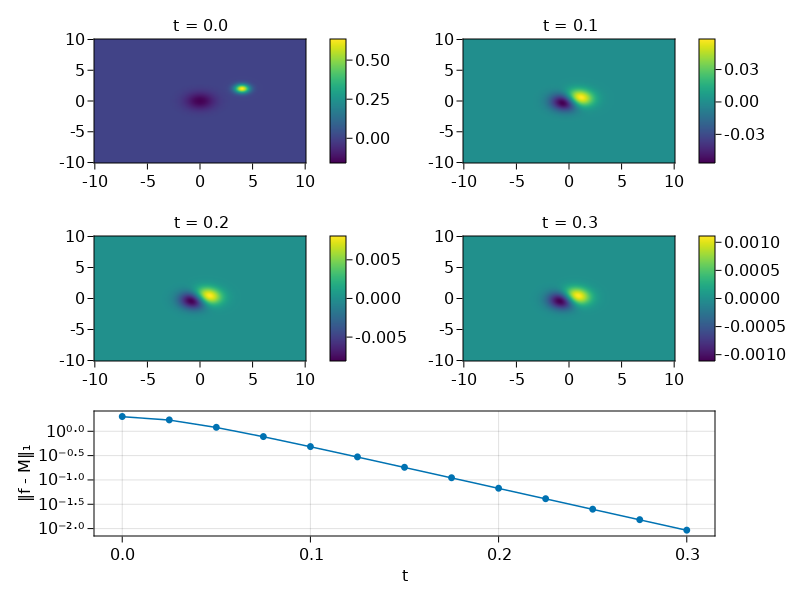}
    \caption{$\epsilon = 0.05$}
  \end{subfigure}
  \begin{subfigure}{0.8\textwidth}
    \centering
    \includegraphics[width=\linewidth]{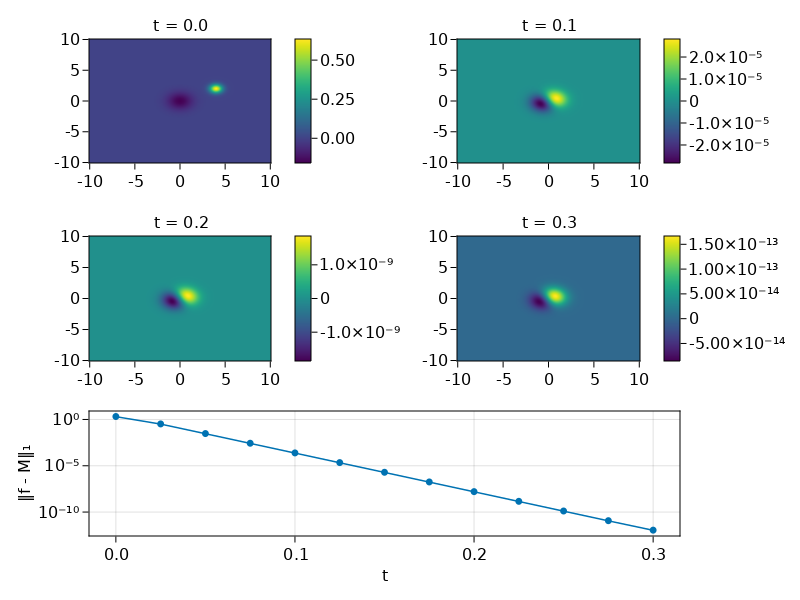}
    \caption{$\epsilon = 0.01$}
  \end{subfigure}
  \caption{Comparison of the characteristic relaxation rates of a cold beam in two velocity dimensions. Heatmaps are the raw difference between $f$ and $f^{0}$, demonstrating convergence to the Maxwellian.}
\end{figure}

%% file: s_step.tex
In this section we motivate the choice of a Forward Euler timestep for the stiff
part of the $S$ flow in (\ref{eqn:S_forward_euler}). Recall that the $S$ step
evolves just the singular values of the solution backwards in time, and comes in
between the $K$ and $L$ steps:
\begin{align}
  \text{$K$ step:} \qquad g^{n} &\xrightarrow{X^{n} \rightarrow X^{n+1},\ S^{n} \rightarrow S^{1}} g^{1} \nonumber \\
  \text{$S$ step:} \qquad g^{1} &\xrightarrow{S^{1} \rightarrow S^{2}} g^{2} \\
  \text{$L$ step:} \qquad g^{2} &\xrightarrow{V^{n} \rightarrow V^{n+1},\ S^{2} \rightarrow S^{n+1}} g^{n+1} \nonumber
\end{align}
Because the $S$ step is backwards in time, strictly speaking it is an ill-posed
ODE. However, empirically we find that an IMEX step for $K$ followed by a
Forwards Euler step for $S$ is stable. To motivate this with a heuristic
argument, consider an initial condition that is uniform in space, with vanishing
electric field and current. We also assume that the solution begins in local
equilibrium.
In terms of our low-rank method, we take
\begin{equation}
g(x, v, t_{0}) = 1, \quad E(x) = J(x) = 0.
\end{equation}
The low-rank decomposition of $g$ gives us $S_{ij} = \delta_{i1}\delta_{j1}$,
i.e. $S_{11}$ is the only nonzero entry of $S$.
With no spatial dependence, all but the collisional term of
(\ref{eqn:partial_t_K_matrix_form}) drop out, and we are left with
\begin{equation}
  \partial_{t} K_{j} = \frac{1}{\epsilon} \sum_{l} d^{1}_{jl} K_{l},
\end{equation}
where $d^{1}_{jl}$ is defined in (\ref{eqn:defn_d_ints}).
The time evolution equation (\ref{eqn:partial_t_S}) for $S$ also simplifies:
\begin{equation}
  \partial_{t} S_{ij} = -\frac{1}{\epsilon} \sum_{kl} \delta_{ik} d^{1}_{jl} S_{kl} = -\frac{1}{\epsilon} \sum_{l} d^{1}_{jl} S_{il}.
\end{equation}

It is useful to rewrite these equations in matrix form.
Define the matrices $\bX\in \mathbb{R}^{N_x\times r}$, $\bS\in \mathbb{R}^{r\times r}$, $\bV\in \mathbb{R}^{N_v\times r}$.
Then $g = \bX \bS \bV^{T} = \bK \bV^{T}$. Further define $\bD = \{ d^{1}_{jl}\}$. The $K$ flow and $S$ flow are given by
\begin{align}
  \partial_{t} \mathbf{K} &= \frac{1}{\epsilon} \bK\bD^T, \\
  \partial_{t} \mathbf{S} &= -\frac{1}{\epsilon} \mathbf{S} \bD^T.
\end{align}

During the $K$ step and the $S$ step, the respective time derivatives of
$g$ are equal and opposite:
\begin{align}
  \text{$K$
  step:} \qquad &\partial_{t}g = (\partial_{t} \mathbf{K}) \mathbf{V}^T = \frac{1}{\epsilon} \bK\bD^T \mathbf{V}^T=\frac{1}{\epsilon} \bX\bS\bD^T \mathbf{V}^T, \\
  \text{$S$
  step:} \qquad &\partial_{t} g = \mathbf{X} (\partial_{t} \mathbf{S})\mathbf{V}^T= -\frac{1}{\epsilon} \mathbf{X} \mathbf{S} \bD^T \mathbf{V}^T.
\end{align}
It follows that at the continuous level, our low-rank approximation has the property that
$g^{n} = g^{2}$ for spatially homogeneous starting point $g^{n}$. We choose our time
discretization to preserve this invariant.
Using the backward Euler for the $K$ flow results in
\begin{equation}
\bK^{n+1}=\bK^n\left( I - \frac{\Delta t}{\epsilon} \bD^T\right)^{-1},
\end{equation}
while using the forward Euler for the $S$ flow results in
\begin{equation}
\bS^{2}=\bS^1\left( I - \frac{\Delta t}{\epsilon} \bD^T\right).
\end{equation}
Therefore,
\begin{align*}
  g^{2} &= \mathbf{X}^{n+1} \mathbf{S}^{2} (\mathbf{V}^{n})^T \\
  &= \mathbf{X}^{n+1} \bS^1\left( I - \frac{\Delta t}{\epsilon} \bD^T\right) (\mathbf{V}^{n})^T \\
        &= \bK^{n+1}\left( I - \frac{\Delta t}{\epsilon} \bD^T\right) (\mathbf{V}^{n})^T \\
&=\bK^n\left( I - \frac{\Delta t}{\epsilon} \bD^T\right)^{-1} \left( I - \frac{\Delta t}{\epsilon} \bD^T\right) (\mathbf{V}^{n})^T \\
  &= g^{n}.
\end{align*}
While we have exact cancellation of the $K$ and $S$ flows for this spatially
homogeneous equilibrium, each of these substeps is quite large when taken
individually. The stiffness of the $\epsilon^{-1}$ term means that it is
quite important to preserve this cancellation at the numerical level, otherwise
the method is unable to hold even a spatially homogeneous equilibrium. For example,
if one uses an IMEX step for the $S$ flow, one finds
\begin{align*}
  g^{2} &= \bK^{n} \left(I - \frac{\Delta t}{\epsilon} \bD^{T} \right)^{-1} \left(I + \frac{\Delta t}{\epsilon} \bD^{T} \right)^{-1} (\bV^{n})^{T} \\
  &= \bK^{n} \left(I - \left( \frac{\Delta t}{\epsilon} \right)^{2} (\bD^{T})^2 \right)^{-1} (\bV^{n})^{T},
\end{align*}
which is very far from the identity indeed.
The preceding argument is not a rigorous justification of our choice of
timestepping scheme for general solutions $g$. We simply wish to highlight one
subtle numerical aspect of the projector-splitting approach which implementors
should be aware of.

\subsection{Comparison to the BGK operator}
In \cite{einkemmerEfficientDynamicalLowRank2021}, the authors successfully used
an IMEX step to advance the $S$ flow. The issues raised above do not arise
for the BGK-type operators considered there, as we demonstrate here with a simple example.
A spatially homogeneous equation with BGK-type collision operator is
\begin{equation*}
  \partial_{t} f = \frac{1}{\epsilon}(M - f),
\end{equation*}
or, using the fact that $M$ is constant for a spatially homogeneous problem,
\begin{equation*}
  \partial_{t} g = \frac{1}{\epsilon}(1 - g).
\end{equation*}
Projecting this onto the low-rank approximation gives the following subflows for $K$ and $S$:
\begin{equation*}
  \partial_{t} K_{j} = \frac{1}{\epsilon}(\langle V_{j} \rangle_{v} - K_{j}), \quad \partial_{t}S_{ij} = -\frac{1}{\epsilon}(\langle X_{i} V_{j} \rangle_{xv} - S_{ij}).
\end{equation*}
An IMEX (backwards Euler) step for each of these subflows will give
\begin{align*}
  \bK^{n+1} &= \left( 1 + \frac{\Delta t}{\epsilon} \right)^{-1} \left( \bK^{n} + \frac{\Delta t}{\epsilon} \langle (\bV^{n})^{T}\rangle_{v} \right), \\
  \bS^{2} &= \left(1 - \frac{\Delta t}{\epsilon} \right)^{-1} \left( \bS^{1} - \frac{\Delta t}{\epsilon} \langle \bX^{n+1}(\bV^{n})^{T} \rangle_{xv} \right).
\end{align*}
Plugging these into the expression for $g^{2}$, we find
\begin{align*}
  g^{2} &= \bX^{n+1} \bS^{2}(\bV^{n})^{T} \\
        &= \bX^{n+1}\left( \bS^{1} - \frac{\Delta t}{\epsilon} \langle \bX^{n+1} (\bV^{n})^T \rangle_{xv} \right) \left( 1 - \frac{\Delta t}{\epsilon} \right)^{-1} (\bV^{n})^{T} \\
        &= \bK^{n+1}\left( 1 - \frac{\Delta t}{\epsilon} \right)^{-1}(\bV^{n})^{T} - \bX^{n+1} \left( \frac{\epsilon}{\Delta t} - 1\right)^{-1} \langle \bX^{n+1} (\bV^{n})^{T} \rangle_{xv} (\bV^{n})^{T} \\
  &= \bX^{n+1} \langle 1, \bX^{n+1} (\bV^{n})^{T} \rangle_{xv} (\bV^{n})^{T} + O\left( \frac{\epsilon}{\Delta t} \right).
\end{align*}
Note that we have dropped terms of order $\epsilon/\Delta t$, to illustrate that $g$ is driven to within $\epsilon/\Delta t$ of its equilibrium value, which is 1 (projected onto the low-rank bases).
Because the BGK operator on $g$ is affine, rather than linear, there is no cancellation, but the IMEX approach
for both flows poses no problems in the $\epsilon \rightarrow 0$ limit.